\documentclass[12pt]{article}
\usepackage{stmaryrd}
\usepackage{mathrsfs}
\usepackage[centertags]{amsmath}
\usepackage{amsfonts, dsfont}
\usepackage{amssymb}
\usepackage{amsthm}
\usepackage{newlfont}
\usepackage{bezier, amsxtra, amsbsy, amsgen, amsopn, amstext}
\usepackage{hyperref}

\input xy
\xyoption{all}

\theoremstyle{plain}
\newtheorem{thm}{Theorem}[subsection]
\newtheorem{cor}[thm]{Corollary}
\newtheorem{lem}[thm]{Lemma}
\newtheorem{prop}[thm]{Proposition}
\newtheorem*{thm2}{Theorem}

\theoremstyle{definition}
\newtheorem{defn}[thm]{Definition}

\newtheorem{remark}[thm]{Remark}

\newtheorem*{ack}{Acknowledgments}

\newcommand{\bd}{\begin{defn}}
\newcommand{\ed}{\end{defn}}
\newcommand{\bl}{\begin{lem}}
\newcommand{\el}{\end{lem}}
\newcommand{\bp}{\begin{prop}}
\newcommand{\ep}{\end{prop}}
\newcommand{\bt}{\begin{thm}}
\newcommand{\et}{\end{thm}}
\newcommand{\bc}{\begin{cor}}
\newcommand{\ec}{\end{cor}}
\newcommand{\br}{\begin{remark}}
\newcommand{\er}{\end{remark}}
\newcommand{\bdi}{\begin{diagram}}
\newcommand{\edi}{\end{diagram}}
\newcommand{\beq}{\begin{eqn}}
\newcommand{\eeq}{\end{eqn}}
\newcommand{\ba}{\begin{array}}
\newcommand{\ea}{\end{array}}
\newcommand{\bpf}{\begin{proof}}
\newcommand{\epf}{\end{proof}}

\newcommand{\Z}{\mathds{Z}}
\newcommand{\Q}{\mathds{Q}}
\newcommand{\Zp}{\mathds{Z}_{p}}
\newcommand{\Qp}{\mathds{Q}_{p}}
\newcommand{\RR}{\mathds{R}}
\newcommand{\al}{\alpha}

\newcommand{\ga}{\gamma}
\newcommand{\Ga}{\Gamma}

\newcommand{\e}{\varepsilon}
\newcommand{\w}{\omega}

\newcommand{\s}{\sigma}
\newcommand{\Si}{\Sigma}
\newcommand{\La}{\Lambda}
\newcommand{\la}{\lambda}
\newcommand{\id}{\mr{id}}

\newcommand{\So}{\Si^{\circ}}
\newcommand{\Om}{\Omega}
\newcommand{\Oo}{\Om^{\circ}}
\newcommand{\Lo}{\La^{\circ}}
\newcommand{\ft}{\mr{ft}}
\DeclareMathOperator{\Mod}{Mod}

\newcommand{\lri}{\stackrel{\sim}{\lra}}

\newcommand{\m}{\mf{m}}
\newcommand{\M}{\mf{M}}
\newcommand{\F}{\mathscr{F}}
\newcommand{\U}{\mathscr{U}}
\newcommand{\C}{\mc{C}}

\newcommand{\II}{\mathscr{I}}
\newcommand{\Sa}{\mc{S}}
\newcommand{\A}{\mf{A}}
\newcommand{\HH}{\mathscr{H}}

\DeclareMathOperator{\Gal}{Gal}
\DeclareMathOperator{\Hom}{Hom}
\DeclareMathOperator{\Ext}{Ext}
\DeclareMathOperator{\RHom}{\mathbf{R}\!\Hom}
\DeclareMathOperator{\Tor}{Tor}
\DeclareMathOperator{\Cone}{Cone}
\DeclareMathOperator{\res}{res}
\DeclareMathOperator{\cTor}{c-Tor}

\newcommand{\ot}{\otimes}
\newcommand{\cotimes}[1]{\,\hat{\ot}_{#1} \,}

\newcommand{\ilim}{\displaystyle \mathop{\varinjlim}\limits}
\newcommand{\plim}{\displaystyle \mathop{\varprojlim}\limits}

\newcommand{\cts}{\mr{cts}}

\newcommand{\DD}{\mathbf{D}}
\newcommand{\R}{\mathbf{R}}
\newcommand{\Li}{\mathbf{L}}

\newcommand{\dotimes}[1]{\, \hat{\ot}{}_{#1}^{\Li} \,}

\newcommand{\ra}{\rightarrow}
\newcommand{\lra}{\longrightarrow}
\newcommand{\tha}{\twoheadrightarrow}

\newcommand{\hra}{\hookrightarrow}

\newcommand{\sbs}{\subseteq}

\newcommand{\mf}{\mathfrak}
\newcommand{\mr}{\mathrm}
\newcommand{\mc}{\mathcal}

\newcommand{\ps}[1]{\llbracket #1 \rrbracket}

\begin{document}

\title{Nekov\'a\v{r} duality over $p$-adic Lie extensions\\ of global fields}
\author{Meng Fai Lim and Romyar T. Sharifi}
\maketitle

\section{Introduction}

\subsection{Duality}

In \cite{Ne}, Nekov\'a\v{r} gave formulations of analogues of Tate and Poitou-Tate duality for finitely
generated modules over a complete commutative local Noetherian ring $R$ with finite residue field of
characteristic a fixed prime $p$.  In the usual formulation of these dualities, one takes the Pontryagin dual,
which does not in general preserve the property of finite generation.  Nekov\'a\v{r}
takes the dual with respect to a dualizing complex of Grothendieck so as to have a duality between bounded complexes of $R$-modules with finitely generated cohomology groups.
This paper is devoted to a generalization of this result to the
setting of nonabelian $p$-adic Lie extensions.

Recall that a dualizing complex $\w_R$ is a bounded complex of $R$-modules
with cohomology finitely generated over $R$ that has
the property that for every complex $M$ of
finitely generated $R$-modules, the Grothendieck dual $\RHom_R(M,\w_R)$ in
the derived category of $R$-modules
$\DD(\Mod_R)$ has finitely generated cohomology,
and moreover, the canonical morphism
\[ M\lra \RHom_R(\RHom_{R}(M,\w_R), \w_R)\] is an
isomorphism in $\DD(\Mod_R)$. Such a complex exists and is unique up
to quasi-isomorphism and translation (see \cite{RD}).

One can choose $\w_R$ to be a bounded complex of injectives, in
which case the derived homomorphism complexes are represented by the
complexes of homomorphisms themselves. If $R$ is regular, then $R$
itself, as a complex concentrated in degree $0$, is a dualizing
complex, but $R$ is not in general $R$-injective.  If $R = \Zp$, for
instance, then the complex $[\Qp \to \Qp/\Zp]$ concentrated in
degrees $0$ and $1$
provides a complex of injective $\Zp$-modules
quasi-isomorphic to $\Zp$.

Let us explain Nekov\'a\v{r}'s theorem, and our generalization of
it, in the setting of Poitou-Tate duality. Let $F$ be a global field
with characteristic not equal to $p$. Let $S$ be a finite set of
primes of $F$ that, if $F$ is a number field, contains all primes
over $p$ and any real places, and let $G_{F,S}$ denote the Galois
group of the maximal unramified outside $S$ extension of $F$.
We remark that a finitely generated $R$-module has a canonical
topology arising from its filtration by powers of the maximal ideal of $R$.
Let $T$ be a bounded complex in the category of finitely generated
(topological) $R$-modules with $R$-linear continuous $G_{F,S}$-actions.
We use $\R\Ga(G_{F,S},T)$ to denote the object
in $\DD(\Mod_R)$ corresponding to the complex of continuous
$G_{F,S}$-cochains with coefficients in $T$, and we use
$\R\Ga_{(c)}(G_{F,S},T)$ to denote the derived object attached to
the complex of continuous compactly supported cochains (using Tate
groups for real places in its definition as a cone: see Section \ref{fgic}).

There exists a bounded complex $T^*$
of finitely generated $R$-modules with $R$-linear continuous $G_{F,S}$-actions that
represents $\RHom_R(T,\w_R)$ in the derived category of such modules.
Nekov\'a\v{r}'s duality theorem is then as follows
(see \cite[Proposition 5.4.3(ii)]{Ne}).\footnote{
Actually, Nekov\'a\v{r}'s longer treatise is particularly concerned with a generalization of this
duality that takes place between Selmer complexes, which we do not address in this article.}

\begin{thm2}[Nekov\'a\v{r}]
    We have an isomorphism
    $$
            \R\Ga(G_{F,S}, T) \lri \RHom_R\!\Big(\R\Ga_{(c)}(G_{F,S}, T^{*}(1)), \w_R\Big)[-3]
    $$
    in the derived category of finitely generated $R$-modules.
\end{thm2}

In this paper, we consider a generalization of this to the setting
of noncommutative Iwasawa theory.  Suppose that $F_{\infty}$ is a
$p$-adic Lie extension of $F$ contained in $F_S$. We denote by $\Ga$
the Galois group of the extension $F_{\infty}/F$, and we let $\La =
R\ps{\Gamma}$ denote the resulting Iwasawa algebra
over $R$.  For a finitely generated $R$-module $T$ with a continuous
$R$-linear $G_{F,S}$-action, we define a finitely generated
$\La$-module $\F_{\Ga}(T)$ with a continuous $\Lambda$-linear
$G_{F,S}$-action by
\[
    \F_{\Ga}(T) = \La^{\iota} \ot_R T,
\]
where the superscript $\iota$ denotes that an element of $\Gamma$ in
$\La$ acts here on $\La$ by right multiplication by its inverse, and
where $G_{F,S}$ acts on $\La$ through the quotient map $G_{F,S} \to
\Gamma$ by left multiplication and then diagonally on the tensor
product.

The $G_{F,S}$-cohomology of $\F_{\Ga}(T)$ is interesting in that a version of Shapiro's lemma
provides natural isomorphisms of continuous cohomology groups
\[
    H^n(G_{F,S}, \F_{\Ga}(T)) \cong \plim_{\al}
    H^n(G_{F_{\alpha},S}, T)
\]
for every $n \ge 0$, where the limit is taken over $\alpha$ indexing
the finite Galois extensions $F_{\al}$ of $F$ that are contained in
$F_{\infty}$ (cf.\ \cite[Lemma 5.3.1]{Lim-adic}). That is, the
cohomology groups of $\F_{\Ga}(T)$ are the Iwasawa cohomology groups
of the module $T$ itself for the extension $F_{\infty}/F$, and this
identification is one of $\La$-modules. We also have the analogous
statements for compactly supported cohomology, as seen in
\cite[Proposition 5.3.3]{Lim-adic}. Therefore, we can reduce the
question of finding dualities among Iwasawa cohomology groups of
compact modules to that of obtaining dualities among cohomology
groups for $G_{F,S}$ itself. Note that, however, it is not even a
priori clear in this setting that $H^n(G_{F,S},\F_{\Ga}(T))$ is a
finitely generated $\La$-module, let alone that we can find such a
duality of $\La$-modules.

The following is our main theorem (cf.\ Theorem \ref{Global duality La}).

\begin{thm2}
We have isomorphisms
\begin{eqnarray*}
    &\R\Ga(G_{F,S}, \F_{\Ga}(T))\lra
    \RHom_{\Lo}\!\Big(\R\Ga_{(c)}(G_{F,S},\F_{\Ga}(T^{*})^{\iota}(1)),\La\ot_R^{\Li}\w_R\Big)[-3]&\\
    &\R\Ga_{(c)}(G_{F,S}, \F_{\Ga}(T))\lra
    \RHom_{\Lo}\!\Big(\R\Ga(G_{F,S},\F_{\Ga}(T^{*})^{\iota}(1)),\La\ot_R^{\Li}\w_R\Big)[-3]&
\end{eqnarray*}
in the derived category of finitely generated $\La$-modules.
\end{thm2}

Nekov\'a\v{r} proves the above theorem
in the case that $\Gamma$ is abelian \cite[Theorem 8.5.6]{Ne}.  In fact, in that case, it is an
almost immediate consequence of his above-mentioned theorem.
That is, suppose for instance that $\Gamma \cong \Zp^r$
for some $r \ge 1$.  Then $\La$ is a complete commutative Noetherian local ring with
finite residue field of characteristic $p$.  Moreover, its dualizing
complex is isomorphic to $\La \ot_R^{\Li} \w_R$ in the
derived category of $\La$-modules \cite[Lemma 8.4.5.6]{Ne}. Therefore, the commutative theory
described above applies to $\La$, and Nekov\'a\v{r} is able to deduce the result from this.
On the other hand, since we are working with nonabelian $\Ga$ and
hence noncommutative $\La$, we do not know that there exists a (nice enough)
dualizing complex, and so the proof of our main theorem takes a different route.  The idea
is a simple one, though the proof is rather involved: after reducing to the case that $R$
is regular and $\Zp$-flat, we perform an inductive argument,
using the grading on $\La$ arising from the powers of its augmentation ideal, to deduce
our result from Nekov\'a\v{r}'s.

We remark that, in their manuscript on the noncommutative main conjecture, Fukaya and Kato stated
an analogue of our main theorem, with $\La$ replacing
$\La \ot_R^{\Li} \w_R$ \cite[(1.6.12)]{FK}, which in turn generalized a result of Burns and Flach
\cite[Lemma 12(b)]{bf} for a narrower class of rings.
The result of Fukaya-Kato applies to a more general class of (adic) rings $\La$ than ours and replaces
$\F_{\Ga}(T)$ and $\F_{\Ga}(T^*)^{\iota}$ by a bounded complex of $\La[G_{F,S}]$-modules $X$
and its $\La$-dual $\Hom_{\La}(X,\La)$.  However, in order to be able to work
with the $\La$-dual, they assume that $X$ consists of (finitely generated) projective
$\La$-modules.
In the case of Iwasawa cohomology that we study, the complex $T$
need not be quasi-isomorphic to a bounded complex of
$R[G_{F,S}]$-modules that are projective and finitely generated over
$R$.  Moreover, if $R$ is Gorenstein, then $R$ serves as an
$R$-dualizing complex, and our result reduces to a duality with
respect to $\La$ itself, as in the result of Fukaya-Kato.  We also note that Vauclair proved a noncommutative duality theorem
for induced modules in the case that $R = \Zp$ and $T$ is
$\Zp$-free, via a rather different method \cite[Theorem 6.4]{vauclair}.

\subsection{An application}

Since applications of our main result are not discussed in the body of this paper,
we end the introduction with an indication of one setting in which our results naturally apply.
Fix $N \ge 1$ not divisible by $p$, and suppose that $p \ge 5$.  Let $\Z_{p,N}$ denote
the inverse limit of the rings $\Z/Np^r\Z$ over $r \ge 1$.  Hida's ordinary cuspidal $\Zp$-Hecke algebra $\mf{h}$ of level $Np^{\infty}$
is a direct product of local rings that is free of finite rank over the subalgebra $\Om = \Zp\ps{1+p\Zp}$ of the algebra $\Zp\ps{\Z_{p,N}^{\times}/\langle -1 \rangle}$
of diamond operators in $\mf{h}$ \cite[Theorem 3.1]{hida}.  Hida showed that
the $\mf{h}$-module $\mc{S}$ of $\Om$-adic cusp forms is $\Om$-dual to $\mf{h}$ (see \cite[Theorem 2.5.3]{ohta} for a proof), from which it follows that $\mc{S}$ is a dualizing complex for $\mf{h}$.

The inverse limit
$$
    \HH = \plim_r H^1_{\text{\'et}}(X_1(Np^r)_{/\bar{\Q}},\Zp(1))^{\text{ord}}
$$
of ordinary parts of \'etale cohomology groups of modular curves is an $\mf{h}[G_{\Q,S}]$-module
for the dual action of $\mf{h}$, where $S$ is the set of primes dividing $Np\infty$.
As an $\mf{h}$-module, Ohta showed in \cite{ohta-ord2} that
$\HH$ is an extension of $\mc{S}$ by $\mf{h}$.
Since $\mf{h}$ is not always Gorenstein (in Eisenstein components: see \cite[Corollary 4.2.13]{ohta-comp} for conditions), it is not at all clear that $\HH$ is quasi-isomorphic to a bounded complex of $\mf{h}[G_{\Q,S}]$-modules that are finitely generated and projective over $\mf{h}$, precluding the use of the duality result of Fukaya-Kato.  On the other hand,
there is a perfect, $\mf{h}$-bilinear pairing $\HH \times \HH \to \mc{S}(1)$
which is $G_{\Q,S}$-equivariant for the action of
$\sigma \in G_{\Q,S}$ on $\mc{S}$ via the diamond operator $\langle \chi(\sigma) \rangle$,
where $\chi \colon G_{\Q,S} \to \Z_{p,N}^{\times}$ is the cyclotomic character (see \cite[Section 1.6]{FK-proof}, where it is quickly derived from a pairing of Ohta's \cite[Definition 4.1.17]{ohta}).

Now fix a $p$-adic Lie extension $F_{\infty}$ of $\Q$ that is unramified outside $S$ and contains $\Q(\mu_{Np^{\infty}})$, and set $\Gamma = \Gal(F_{\infty}/\Q)$.
The complex $\R\Ga_{(c)}(G_{\Q,S},\F_{\Ga}(\HH))$ is the subject of the noncommutative
Tamagawa number conjecture of Fukaya and Kato for $\HH$
(though to be precise, said conjecture is only formulated in the case that $\HH$ is $\mf{h}$-projective), which
is directly related to the noncommutative main conjecture for ordinary $\Omega$-adic
cusp forms (see \cite[Sections 2.3 and 4.2]{FK}).
It is also perhaps worth remarking that, in the commutative setting, the first Iwasawa cohomology group $H^1(G_{\Q,S},\F_{\Ga}(\HH))$ contains zeta elements constructed out of Kato's Euler system (see \cite[Section 12]{kato} or \cite[Section 3.2]{FK-proof}).

For a finitely generated $\La = \mf{h}\ps{\Ga}$-module $A$, let us use $A\langle \chi \rangle$ to denote the $\La$-module that is $A$ as an $\mf{h}$-module but for which the original $\mf{h}$-linear action of $\gamma \in \Gamma$ on $A$ has been twisted by multiplication by the diamond operator $\langle \chi(\tilde{\gamma}) \rangle$,
for any lift $\tilde{\gamma} \in G_{\Q,S}$ of $\gamma$.
Using the pairing on $\HH$, our main result can be seen to yield two interesting isomorphisms in the derived category of finitely generated modules over, including
$$
    \R\Ga_{(c)}(G_{\Q,S},\F_{\Ga}(\HH)) \lri
    \RHom_{\Lo}(\R\Ga(G_{\Q,S},\F_{\Ga}(\HH)^{\iota}),\La \otimes_{\mf{h}}^{\Li}
    \mc{S})\langle \chi \rangle[-3].
$$
Perhaps more concretely, we have a spectral sequence
$$
    \Ext^i_{\Lo}\Big(H^{3-j}(G_{\Q,S},\F_{\Ga}(\HH)^{\iota}),
    \Lambda \otimes_{\mf{h}} \mc{S}\Big)\langle \chi \rangle
    \Longrightarrow H^{i+j}_{(c)}(G_{\Q,S},\F_{\Ga}(\HH)).
$$

We remark that one has a similar result with $X_1(Np^r)$ replaced by $Y_1(Np^r)$: in this case, the pairing on $\HH$ to $\mc{S}$ is replaced by a pairing between ordinary parts
of cohomology and compactly supported cohomology groups to ordinary $\Omega$-adic modular forms.
We also note that Fouquet has constructed an analogue of Ohta's pairing
for towers of Shimura curves attached to indefinite quaternion algebras over totally real fields \cite[Proposition 2.8]{fouquet}, providing a related setting for an application of our main
result.

\begin{ack}
    The material presented in this article forms a revised and extended version of the part of
    the Ph.D.\ thesis \cite{Lim-thesis} of the first author that was a collaboration between the
    two authors.  The results of the rest of the thesis can be
    found in the article \cite{Lim-adic}.
    The first author would like to thank Manfred Kolster for his encouragement.
    The second author would like to thank Takako Fukaya, Uwe Jannsen, Jan Nekov\'a\v{r},
    Otmar Venjakob, and Amnon Yekutieli for helpful discussions.
    He was supported by NSF Grant DMS-0901526, and his work was partially conducted
    during visits to the Newton Institute in 2009 and IHES in 2010.
    The authors also thank the anonymous referee for pointing out a mistake in an earlier version
    of this work.
\end{ack}

\section{Preliminaries}

In this section, we lay out some facts regarding modules over
profinite rings that we shall require later in the paper. Elsewhere
in the paper, we use some of the standard terminology and
conventions regarding complexes, cones, shifts, total complexes of
homomorphisms and tensor products, and derived categories, as can be
found in \cite{Ne}, or in \cite{Lim-adic} off of which this paper
builds.  In particular, all sign conventions will be as in the
latter two papers.  For further background on derived categories, we
suggest the excellent book \cite{Wei}, as well as the more advanced text \cite{KS}.

\subsection{Tensor products and homomorphism groups}
\label{tens hom}

In this subsection, $R$ is a commutative ring.  By an $R$-algebra,
we will mean a ring with a given homomorphism from $R$ to its
center.  In this section, we construct derived bifunctors of
homomorphism groups and tensor products for
bimodules over $R$-algebras and study a few isomorphisms
that result from these constructions.

We use $\Lo$ to denote the opposite ring to an $R$-algebra $\La$.
Given two $R$-algebras
$\La$ and $\Si$, we are interested in a
subclass of the class of $\La$-$\Si$-bimodules, namely, the class of
$\La$-$\Si$-bimodules with the extra property that the left
$R$-action coincides with the right $R$-action.
We can (and shall)
identify the category of such $\La$-$\Si$-bimodules with the
category of $\La\ot_{R}\So$-modules, and there are natural exact
forgetful functors
\begin{eqnarray*}
    \res_{\La} \colon \Mod_{\La\ot_{R}\So}\lra \Mod_{\La}
    &\mr{and}&
    \res_{\So} \colon \Mod_{\La\ot_{R}\So}\lra
    \Mod_{\So},
\end{eqnarray*}
which extend to exact functors on the derived
categories.
One observes that the categories $\Mod_{\La\ot_{R}R^{\circ}}$ and
$\Mod_{\La}$ are equivalent, as are
$\Mod_{R\ot_{R}\So}$ and $\Mod_{\So}$.
We have following lemma, which extends that of \cite[Lemma 2.1]{Ye}.

\bl \label{abstract flat}\ \begin{enumerate} \item[$(a)$] If $\Si$
is a projective $($resp.,\ flat$)$ $R$-algebra, then $\res_{\La}$
preserves projective $($resp.,\ flat$)$ modules.
Similarly, if $\La$ is a projective $($resp.,\ flat$)$ $R$-algebra, then
$\res_{\So}$ preserves projective $($resp.,\ flat$)$ modules.

\item[$(b)$]  If $\Si$ is a flat $R$-algebra, then $\res_{\La}$
preserves injectives.  Similarly, if $\La$ is a flat $R$-algebra, then
$\res_{\So}$ preserves injectives.
\end{enumerate}\el

\bpf (a) We prove the statements of the first sentence, those of the second being
a consequence.  Suppose that $\Si$ is a projective $R$-algebra. Since
projective modules are the summands of free modules, it
suffices to show that $\La\ot_{R}\So$ is a projective
$\La$-module. Since $\Si$ is a central $R$-algebra, we have
$\Si\cong \So$ as $R$-modules. Therefore, we have
$\So \oplus P\cong L$ for some projective $R$-module $P$ and
free $R$-module $L$. Then $\La\ot_{R}\So$ is a direct
summand of $\La\ot_{R}L$, which is a free $\La$-module. Hence,
$\La\ot_{R}\So$ is a projective $\La$-module.

Now suppose that $\Si$ is flat over $R$. Since flat modules are direct
limits of finitely generated free modules (see \cite[Theorem 4.34]{Lam})
and tensor products preserve direct limits, it suffices
to show that $\La\ot_R \So$ is a flat $\La$-algebra. Since
$\Si$ is flat over $R$, we have that $\Si$ is a direct limit of
finitely generated free $R$-modules, which implies that $\La\ot_R
\So$ is a direct limit of finitely generated free
$\La$-modules.

 \vspace{2ex}
(b) We shall prove this for $\res_{\La}$, the case of $\res_{\So}$
being a consequence. The functor from $\Mod_{\La}$ to $\Mod_{\La \ot_R \So}$
sending $M$ to $M \ot_R \So$ is exact by our assumption and is
left adjoint to the functor $\res_{\La}$. The conclusion then follows
from \cite[Prop.\ 2.3.10]{Wei}. \epf

Let us briefly introduce the notions of q-projective
and q-injective complexes of $\La$-modules and several facts regarding them arising
from the work of Spaltenstein,
as can be found in \cite{Sp}, \cite{Kel}, \cite[Section 2.3]{Lip}, and \cite[Chapter 14]{KS}.
A complex $P$ of $\La$-modules is called q-projective
if for every map $g \colon P \to B$ and quasi-isomorphism $s \colon A \to B$
of complexes of $\La$-modules, there exists a map $f \colon P \to A$ of complexes of $\La$-modules
such that $g$ and $s \circ f$ are homotopy equivalent.
We remark that $P$ is q-projective if and only if it is
homotopy equivalent to a direct limit of bounded above complexes
$P_n$ of projective $\La$-modules via maps $\iota_n \colon
P_n \to P_{n+1}$ that are split injective in each degree with quotients
$P_{n+1}/\iota_n(P_n)$ having zero differentials.
In particular, bounded above complexes of projectives are q-projective.
If $A$ is a complex of $\La$-modules, then there exists a quasi-isomorphism
$P \to A$ with $P$ q-projective.

We also have the dual notion of q-injective complexes.   A complex
$I$ of $\La$-modules is q-injective if for every map $f \colon A \to
I$ and quasi-isomorphism $s \colon A \to B$ of complexes of
$\La$-modules, there exists a map $g \colon B \to I$ such that $f$
and $g \circ s$ are homotopy equivalent. A complex $I$ of
$\La$-modules is q-injective if and only if it is homotopy
equivalent to an inverse limit of bounded below complexes $I_n$ of
injective $\La$-modules via maps $\pi_n \colon I_{n+1} \to I_n$ that
are split surjective in each degree with kernels $\ker \pi_n$ having
zero differentials. In particular, bounded below complexes of
injective $\La$-modules are q-injective. If $A$ is a complex of
$\La$-modules, then there exists a quasi-isomorphism $A \to I$ with
$I$ q-injective.

In addition to $\La$ and $\Si$, we now let $\Om$ be an
$R$-algebra. If $A$ is a $\La\ot_R\Oo$-module and $B$
is a $\La\ot_R\So$-module, we give $\Hom_{\La}(A,B)$
the structure of an $\Om \ot_R \So$-module via the
left $\Om$ and right $\Si$-actions
\begin{eqnarray*}
    (\w \cdot f)(a) = f(a\w) &\mr{and}& (f\cdot\s)(a) = f(a)\s
\end{eqnarray*}
for $f \in \Hom_{\La}(A,B)$, $a \in A$, $\w \in \Om$, and $\s \in \Si$.

Moreover, if $A$ is a complex of $\La\ot_R\Oo$-modules and
$B$ is a complex of $\La\ot_R\So$-modules, we define a
complex $\Hom_{\La}(A,B)$ of
$\Om\ot_R\So$-modules by
\[
    \Hom_{\La}^{n}(A,B) = \prod_{i\in\Z}\Hom_{\La}(A^{i}, B^{i+n}),
\]
with the usual differentials, as in \cite[Section 2]{Lim-adic}.

\bp \label{derived Hom}
    Let $A$ be a complex of $\La \ot_R \Oo$-modules, and let $B$ be a
    complex of $\La \ot_R \So$-modules.
    If $\Om$ is a projective $R$-algebra or $\Si$ is a flat $R$-algebra,
    then we have a derived bifunctor
    $$
        \RHom_{\La}(-,-) \colon \DD(\Mod_{\La\ot_{R}\Oo})^{\circ}
        \times \DD(\Mod_{\La\ot_{R}\So}) \lra \DD(\Mod_{\Om \ot_R \So}),
    $$
    and $\Hom_{\La}(A,B)$ represents $\RHom_{\La}(A,B)$ if
    $A$ is q-projective as a complex of $\La$-modules and $\Om$ is $R$-projective
    or $B$ is q-injective as a complex of $\La$-modules and $\Si$ is $R$-flat.
\ep

\bpf
    Let us first assume that $\Om$ is a projective $R$-algebra.
    Since every complex of $\La \ot_R \Oo$-modules is quasi-isomorphic
    to a q-projective complex of $\La \ot_R \Oo$-modules,
    we have a derived functor
    \[ \RHom_{\La}(-,B) \colon  \DD(\Mod_{\La\ot_{R}\Oo})^{\circ}
    \lra \DD(\Mod_{\Om \ot_R \So}) \]
    (see \cite[Corollary 2.3.2.3]{Lip}).
    Suppose that $f \colon B \to B'$ is a
    quasi-isomorphism of complexes of $\La\ot_R\So$-modules.
    Let $\e \colon P \to A$
    be a quasi-isomorphism of complexes of $\La \ot_R \Oo$-modules,
    where $P$ is q-projective.
    Then $\RHom_{\La}(P,B)$ (resp., $\RHom_{\La}(P,B')$) is
    represented by $\Hom_{\La}(P,B)$ (resp.,
    $\Hom_{\La}(P,B')$).  Applying Lemma \ref{abstract flat}(a), we see that $P$ is also a
    q-projective complex of $\La$-modules, so by \cite[Proposition 2.3.8]{Lip},
    the induced map
     \[f_* \colon \Hom_{\La}(P,B) \lra \Hom_{\La}(P,B')\]
     is a quasi-isomorphism of complexes of abelian groups.
     Since $f_*$ is a morphism
     of complexes of $\Om\ot_R\So$-modules, it is a quasi-isomorphism of
     such complexes.  Hence, $f$ induces isomorphisms
     $$
        \RHom_{\La}(A,B) \lra \RHom_{\La}(A,B'),
    $$
    proving the existence of the derived bifunctor.
    Moreover, if $A$ is q-projective as a complex of
    $\La$-modules, then
     \[
        \e_* \colon \Hom_{\La}(A,B) \lra \Hom_{\La}(P,B)
    \]
    is a
    quasi-isomorphism of complexes of $\So$-modules, hence of
    $\Om \ot_R \So$-modules as well.  Thus, $\Hom_{\La}(A,B)$
    represents $\RHom_{\La}(A,B)$, as desired.

    If $\Si$ is $R$-flat,
    the argument for the existence of the derived bifunctor and its computation by $\Hom_{\La}(A,B)$
    in the case that $B$ is q-injective as a complex of $\La$-modules
    is the direct analogue of the above argument, employing
    Lemma \ref{abstract flat}(b).  To see that the resulting derived functor
    agrees with that constructed
    above in the case that both $\Om$ is $R$-projective
    and $\Si$ is $R$-flat, consider a quasi-isomorphism
    $\e \colon P \to A$ as above and a quasi-isomorphism
    $\delta \colon B \to I$ with $I$ a q-injective complex of $\La\ot_R\So$-modules.
    As we have quasi-isomorphisms
    \[
        \SelectTips{cm}{} \xymatrix{
        \Hom_{\La}(P,B) \ar[r]^-{\delta_*} & \Hom_{\La}(P,I)& \ar[l]_-{\e_*}
        \Hom_{\La}(A,I) }
    \]
    of complexes of $\Om \ot_R \So$-modules, the derived functors coincide.
\epf

Note that one always has a map $\Hom_{\La}(A,B) \to
\RHom_{\La}(A,B)$ (canonical up to isomorphism in the derived
category), induced either by a quasi-isomorphism $P \to A$ with $P$ q-projective or a quasi-isomorphism
$B \to I$ with $I$ q-injective.

\bc \label{derived Hom res}
 If $\Om$ is a projective $R$-algebra,
then we have a commutative diagram
 \[
    \SelectTips{cm}{} \xymatrix @C=1in{
     \DD(\Mod_{\La\ot_{R}\Oo})^{\circ}
    \times \DD(\Mod_{\La\ot_{R}\So}) \ar[d]_{\res_{\La}\times \id}
     \ar[r]^-{\RHom_{\La}(-,-)} & \DD(\Mod_{\Om\ot_R\So})
     \ar[d]^{\res_{\So}} \\
   \DD(\Mod_{\La})^{\circ} \times \DD(\Mod_{\La\ot_R\So}) \ar[r]^-{\RHom_{\La}(-,-)} & \DD(\Mod_{\So}).  }
\]
In the case that $\Si$ is a flat $R$-algebra, we have a commutative
diagram
 \[
    \SelectTips{cm}{} \xymatrix @C=1in{
     \DD(\Mod_{\La\ot_{R}\Oo})^{\circ}
    \times \DD(\Mod_{\La\ot_{R}\So}) \ar[d]_{\id\times \res_{\La}}
     \ar[r]^-{\RHom_{\La}(-,-)} & \DD(\Mod_{\Om\ot_R\So})
     \ar[d]^{\res_{\Om}} \\
   \DD(\Mod_{\La\ot_R\Oo})^{\circ} \times \DD(\Mod_{\La}) \ar[r]^-{\RHom_{\La}(-,-)} & \DD(\Mod_{\Om}).
   }
\] \ec

If $A$ is a complex of $\Om\ot_R\Lo$-modules and $B$ is a complex
of $\La \ot_R \So$-modules, we define a complex
$A\ot_{\La}B$ of $\Om \ot_R \So$-modules by
\[ (A\ot_{\La}B)^{n} = \bigoplus_{i\in\Z}A^{i}\ot_{\La}B^{n-i}, \]
again with the usual differentials, as in \cite[Section
2]{Lim-adic}.

We also have a notion of a q-flat complex of $\Lo$-modules (see
\cite[Section 2.5]{Lip} and \cite[Section 5]{Sp} in the case of
commutative rings, the proofs being identical); that is, a
complex of $\Lo$-modules $A$ is said to be q-flat if for every quasi-isomorphism
$B \to C$ of complexes of $\La$-modules, the resulting map $A \ot_{\La} B \to A
\ot_{\La} C$ is also a quasi-isomorphism. In particular,
q-projective complexes of $\Lo$-modules are q-flat, any bounded
above complex of flat $\Lo$-modules is q-flat, and any filtered
direct limit of q-flat complexes is q-flat. As in the case of
homomorphism complexes, the total tensor product induces derived
bifunctors as follows.  We omit the analogous proof.

\bp \label{derived tensor}
    Let $A$ and $B$ be complexes of $\Om \ot_R \Lo$-modules
    and $\La \ot_R \So$-modules, respectively.
    If either $\Om$ or $\Si$ is a flat
    $R$-algebra, then we have a derived bifunctor
    $$
       -\ot^{\Li}_{\La}- \colon
        \DD(\Mod_{\Om\ot_R\Lo})\times\DD(\Mod_{\La \ot_R \So})\lra
        \DD(\Mod_{\Om \ot_R \So}),
    $$
    and $A \ot_{\La} B$ represents $A\ot_{\La}^{\Li}B$ if $A$ is
    q-flat as a complex of $\Lo$-modules and $\Om$ is $R$-flat or if
    $B$ is q-flat as a complex of $\La$-modules and $\Si$ is $R$-flat.
\ep

Proposition \ref{derived tensor} has the following direct corollary.

\bc \label{tens comm}
    If $\Om$ is a flat $R$-algebra, then we have a
    commutative diagram
     \[
        \SelectTips{cm}{} \xymatrix @C=1in{
        \DD(\Mod_{\Om\ot_R\Lo}) \times \DD(\Mod_{\La\ot_R\So}) \ar[d]_{\res_{\Lo}\times\id}
         \ar[r]^-{-\, \ot_{\La}^{\Li}\, -} & \DD(\Mod_{\Om\ot_R\So})
         \ar[d]^{\res_{\So}} \\
       \DD(\Mod_{\La^{\circ}}) \times \DD(\Mod_{\La\ot_R\So}) \ar[r]^-{-\, \ot_{\La}^{\Li}\, -} &
       \DD(\Mod_{\So}),  }
    \]
    and if $\Si$ is a flat $R$-algebra, we have a
    commutative diagram
     \[
        \SelectTips{cm}{} \xymatrix @C=1in{
        \DD(\Mod_{\Om\ot_R\Lo}) \times \DD(\Mod_{\La\ot_R\So}) \ar[d]_{\id\times\res_{\La}}
         \ar[r]^-{-\, \ot_{\La}^{\Li}\, -} & \DD(\Mod_{\Om\ot_R\So})
         \ar[d]^{\res_{\Om}} \\
       \DD(\Mod_{\Om\ot_R\Lo}) \times \DD(\Mod_{\La}) \ar[r]^-{-\, \ot_{\La}^{\Li}\, -} &
       \DD(\Mod_{\Om}).  }
    \]
\ec

We end this section with some general lemmas regarding the passing
of tensor products through homomorphism groups and the resulting
isomorphisms in the derived categories.
Let us use $\bar{m}$ to denote the degree of
an element $m$ of a term of a complex $M$ of modules over a ring.
We fix a fourth $R$-algebra $\Xi$.

\bl \label{derived adj}
    Suppose that $\Xi$ is a flat $R$-algebra and $\Si$ is a
    projective $R$-algebra.  Let $A$ be a
    complex of $\Om \ot_R \Lo$-modules,
    let $B$ be a complex of $\La \ot_R \So$-modules, and
    let $C$ be a complex of $\Om \otimes_R \Xi^{\circ}$-modules.
    Fix a quasi-isomorphism $Q \to A$ of complexes of $\Om \ot_R \Lo$-modules with
    $Q$ q-flat over $\Lo$, a quasi-isomorphism $P \to B$ of complexes of
    $\La \ot_R \So$-modules with $P$ q-projective over $\La$,
    and a quasi-isomorphism
    $C \to I$ of complexes of
    $\Om \ot_R \Xi^{\circ}$-modules with $I$ q-injective over $\Om$.
    Then the adjunction isomorphism
    \begin{eqnarray*}
        &\Hom_{\Om}(A \ot_{\La} P,I) \lra \Hom_{\La}(P,\Hom_{\Om}(A,I))&\\
        &f \mapsto (p \mapsto (-1)^{\bar{a}\bar{p}}f(a \otimes p))&
    \end{eqnarray*}
    induces an isomorphism
    \[
        \RHom_{\Om}(A\ot^{\Li}_{\La}B,C) \lri
        \RHom_{\La}(B, \RHom_{\Om}(A,C))
    \]
    in $\DD(\Mod_{\Si \ot_R \Xi^{\circ}})$,
    as does the adjunction isomorphism
    $$
        \Hom_{\Om}(Q \ot_{\La} B,I) \lra \Hom_{\La}(B,\Hom_{\Om}(Q,I)).
    $$
\el

\bpf
    By Propositions \ref{derived Hom} and \ref{derived tensor},
    $
         \Hom_{\Om}(A \ot_{\La} P,I)
    $
    represents
    $
        \RHom_{\Om}(A\ot^{\Li}_{\La}B,C),
    $
    and
    $
        \Hom_{\La}(P, \Hom_{\Om}(A,I))
    $
    represents
    $
        \RHom_{\La}(B, \RHom_{\Om}(A,C)).
    $
    Therefore, we are reduced for the first part
    to showing that the adjunction map
    is an isomorphism of complexes, and this is standard
    (see \cite[Lemma 2.2]{Lim-adic}).
    The second part follows easily from the first part and the commutative diagram
    $$
        \SelectTips{cm}{} \xymatrix{
        \Hom_{\Om}(Q \ot_{\La} B, I) \ar[r]^-{\sim} \ar[d] & \Hom_{\La}(B,\Hom_{\Om}(Q,I)) \ar[d] \\
        \Hom_{\Om}(Q \ot_{\La} P, I) \ar[r]^-{\sim} & \Hom_{\La}(P,\Hom_{\Om}(Q,I)),
    }
    $$
    in that the left-hand vertical map is a quasi-isomorphism.
\epf

We consider the next isomorphism first on the level of complexes.

\bl \label{tech 2}
        	Let $A$ be a bounded complex of $\Om \ot_R \So$-modules
    	that are flat as $\So$-modules, let $B$ be a complex
    	of $\Xi \ot_R \Lo$-modules that are finitely presented as $\Lo$-modules, and let $C$ be
        	a complex of $\Si \ot_{R} \Lo$-modules.
    	Suppose also that either the terms of $A$ are finitely presented as $\So$-modules
    	or at least one of $B$ and $C$ is bounded above and at least one is bounded below.
        	Then the map
        	\begin{eqnarray*}
                &A \ot_{\Si} \Hom_{\Lo}(B,C) \lra \Hom_{\Lo}(B,A \ot_{\Si} C)& \\
                &a \ot f \mapsto \big(b\mapsto a \ot f(b)\big)&
        	\end{eqnarray*}
    	is an isomorphism of complexes of $\Om \ot_R \Xi^{\circ}$-modules.
\el

\begin{proof}
    That the stated map is a map of complexes is an easy check of the actions and
    signs.
    To see that it is an
    isomorphism, consider first the case that $A$, $B$, and $C$ are concentrated in degree $0$
    and take a presentation
    $$
        F_1 \lra F_2 \lra B \lra 0
    $$
    with $F_1$ and $F_2$ finitely generated free $\Lo$-modules.
    Then the two rightmost vertical arrows in the commutative diagram with exact
    (noting the $\So$-flatness of $A$) rows
    $$
        \SelectTips{cm}{} \xymatrix @R=0.4in @C=0.2in{
            0\ar[r] & A \ot_{\Si} \Hom_{\Lo}(B,C) \ar[d]_{} \ar[r]^{}  &
        A \ot_{\Si}\Hom_{\Lo}(F_2,C) \ar[d]_{} \ar[r]^{} &
        A \ot_{\Si} \Hom_{\Lo}(F_1,C) \ar[d]^{}   \\
            0\ar[r]  &\Hom_{\Lo}(B,A\ot_{\Si}C)  \ar[r]^{}  &
        \Hom_{\Lo}(F_2, A \ot_{\Si} C)
                \ar[r]^{} & \Hom_{\Lo}(F_1, A \ot_{\Si} C) }
    $$
    are isomorphisms of $\Om$-modules
    by universal property of the direct sum and the commutativity of direct sums and tensor
    products, so the other is as well.
    Since the latter map is a morphism of complexes of $\Om \ot_R \Xi^{\circ}$-modules,
    we have the result in the case of modules.

    In the general setting, we have that
    $$
        \Hom_{\Lo}^n(B, A \ot_{\Si} C) =
        \prod_{j \in \Z} \Hom_{\Lo} \!\Big( B^j, \bigoplus_{i \in \Z} A^i \ot_{\Si} C^{n-i+j} \Big).
    $$
    Since $B^j$ is finitely presented over $\Lo$, the argument
    given above
    implies that the latter group is naturally isomorphic to
    $$
        \prod_{j \in \Z} \bigoplus_{i \in \Z}  A^i \ot_{\Si} \Hom_{\Lo}(B^j,C^{n-i+j}).
    $$
    On the other hand, we have
    $$
        (A \ot_{\Si} \Hom_{\Lo}(B,C))^n = \bigoplus_{i \in \Z} A^i \ot_{\Si} \!
        \Big( \prod_{j \in \Z} \Hom_{\Lo}(B^j,C^{n-i+j}) \Big).
    $$
    If the terms of $A$ are finitely presented over $\So$, then we
    can use a finite presentation
    of $A^i$ to see that the latter term is naturally isomorphic to
    $$
        \bigoplus_{i \in \Z} \prod_{j \in \Z} A^i \ot_{\Si} \Hom_{\Lo}(B^j,C^{n-i+j}),
    $$
    and the boundedness of $A$ allows us to commute the direct product and direct sum.
    Finally, if at least one of $B$ and $C$ is bounded above and at least one is bounded below, then
    the products over $j \in \Z$ involve only finitely many nonzero terms
    (again using that $A$ is bounded),
    and the result follows as direct sums commute
    with tensor products and each other.
\end{proof}

In the derived category, we have the following generalization of \cite[Proposition 6.1]{Ven02}.

\bl \label{tech 2 derived}
    Suppose that $\Om$ is a flat $R$-algebra and $\Xi$ is a projective $R$-algebra.
    Let $A$ be a bounded above complex of $\Om \ot_R \So$-modules, let
    $B$ be a complex of $\Xi \ot_R \Lo$-modules, and let $C$ be
    a bounded below complex of $\Si \ot_{R} \Lo$-modules.  Suppose that there exists a
    quasi-isomorphism $Q \to A$ of $\Om \ot_R \So$-modules, where $Q$ is bounded and has terms
    that are flat as $\So$-modules. Suppose also that there is a quasi-isomorphism
    $P \to B$ of complexes of $\Xi \ot_R \Lo$-modules, where $P$ is q-projective as a
    complex of $\Lo$-modules with terms that are finitely presented over $\Lo$.
    We assume also that either $C$ is bounded, $P$ is bounded above, or the terms of
    $Q$
     are finitely presented as $\So$-modules.
    Then the map
        $$
                Q \ot_{\Si} \Hom_{\Lo}(P,C) \lra \Hom_{\Lo}(P,Q \ot_{\Si} C)
    $$
    gives rise to an isomorphism
    \[
        A \ot_{\Si}^{\Li} \RHom_{\Lo}(B,C) \lri
        \RHom_{\Lo}(B, A\ot^{\Li}_{\Si}C)
    \]
    in $\DD(\Mod_{\Om \ot_R \Xi^{\circ}})$.
\el

\begin{proof}
    By Propositions \ref{derived Hom} and \ref{derived tensor},
    $
        Q \ot_{\Si} \Hom_{\Lo}(P,C)
    $
    represents
    $
        A \ot_{\Si}^{\Li} \RHom_{\Lo}(B,C)
    $
    and
    $
        \Hom_{\Lo}(P,Q \ot_{\Si} C)
    $
    represents
    $
         \RHom_{\Lo}(B, A\ot^{\Li}_{\Si}C).
    $
    That the stated map is an isomorphism follows
    from Lemma \ref{tech 2}.
\end{proof}

\begin{remark}
    In this article, we need only compute derived tensor products in the case
    that at least one of the complexes is bounded above and derived homomorphism
    complexes in the case that the target complex is bounded below.
    However, in Section \ref{global duality}, we will nevertheless be forced to
    represent derived homomorphism complexes using q-projective resolutions
    in an unbounded first variable.  Hence, we have given a general treatment.
\end{remark}

\subsection{Modules over group rings}
\label{mod gp rings}

In this subsection, we preserve the notation of
Section \ref{tens hom}, and in addition, let $G$ be a group.  We begin by extending
the notions of derived bifunctors of homorphism groups and
tensor products to incorporate an additional action of $G$.

If $A$ is a complex of $\La[G]$-modules and $B$ is a complex of $(\La\ot_R\So)[G]$-modules,
we give $\Hom_{\La}(A,B)$  the structure of a complex of $\So[G]$-modules via the standard $G$-action
$$(g \cdot f)(a) = g \cdot f(g^{-1}a)$$
for $g \in G$, $a \in A$, and $f \in \Hom_{\La}(A,B)$.  The following is then a slight weakening of
the natural analogue of Proposition \ref{derived Hom}, written in condensed form.

\bp \label{derived Hom G}
    There is a derived bifunctor
    $$
        \RHom_{\La}(-,-) \colon \DD(\Mod_{\La[G]})^{\circ}
            \times \DD(\Mod_{\La\ot_{R}\So[G]}) \lra \DD(\Mod_{\So[G]}).
    $$
    Moreover, $\RHom_{\La}(A,B)$ can be represented by $\Hom_{\La}(A,B)$
    if $A$ is q-projective as a complex of $\La$-modules or, if $\Si$ is $R$-flat,
    $B$ is q-injective as a complex of $\La$-modules.
\ep

\bpf
    Proposition \ref{derived Hom} (with
    $\Om$ replaced by $R[G]^{\circ}$ and $\Si$ replaced by $\So[G]$)
    implies the existence of derived bifunctors
    with the desired properties that take values in the category $\DD(\Mod_{R[G]^{\circ} \ot_R \So[G]})$.
    There is a natural exact functor
    $$
        \Mod_{R[G]^{\circ} \ot_R \So[G]} \lra \Mod_{\So[G]}
    $$
    taking $M$ to the same $\So$-module with new $G$-action
    $$
        g \cdot m = (g^{-1} \ot g) \cdot m
    $$
    for $g \in G$ and $m \in M$.  This in turn induces a functor on derived categories, and
    composition of the above derived functor with this functor yields the result.
\epf

Next, if $A$ is a complex of $(\Om\ot_R\Lo)[G]$-modules and $B$ is a complex
of $(\La \ot_R \So)[G]$-modules, we give
$A\ot_{\La}B$ the structure of a complex of $(\Om \ot_R \So)[G]$-modules via the $G$-action
defined by $g(a \ot b) = ga \ot gb$ for $a \in A$, $b \in B$, and $g \in G$.  We state
a very slight weakening (for brevity) of the analogue of Proposition \ref{derived tensor}.

\bp \label{derived tensor G}
    If $\Om$ or $\Si$ is a flat $R$-algebra, then we have a derived bifunctor
    $$
        -\ot^{\Li}_{\La}- \colon
            \DD(\Mod_{\Om\ot_R\Lo[G]})\times \DD(\Mod_{\La \ot_R \So[G]})\lra
            \DD(\Mod_{\Om \ot_R \So[G]}).
    $$
    Moreover, $A \ot_{\La} B$ represents $A\ot_{\La}^{\Li}B$ if either $\Om$ is $R$-flat and
    A is a q-flat as a complex of $\Lo$-modules or $\Si$ is $R$-flat and
    $B$ is q-flat as a complex of $\La$-modules.
\ep

\begin{proof}
    Proposition \ref{derived tensor} with $\Om$ replaced by $\Om[G]$ and
    $\So$ replaced by $\So[G]$ yields the existence of a derived bifunctor
    $$
        -\ot^{\Li}_{\La}- \colon
            \DD(\Mod_{\Om\ot_R\Lo[G]})\times \DD(\Mod_{\La \ot_R \So[G]})\lra
            \DD(\Mod_{\Om[G] \ot_R \So[G]})
    $$
    which has the desired properties.
    There is a natural exact functor
    $$
        \Mod_{\Om[G] \ot_R \So[G]}
        \lra \Mod_{(\Om \ot_R \So)[G]}
    $$
    that takes an module $M$ to the same $\Om \ot_R \So$-module with new $G$-action
    $$
        g \cdot m = (g \ot g) \cdot m
    $$
    for $g \in G$ and $m \in M$,
    and the desired bifunctor is the resulting composition of derived functors.
\end{proof}

Suppose now that $\La$ is an $R$-algebra and that we are given a homomorphism
$\chi \colon G \to \mr{Aut}_{\Lo}(\La)$, allowing us to endow $\La$ with the structure of a
$\Lo[G]$-module.  We denote the resulting module by ${}_{\chi}\La$.  (If $\chi$ is trivial,
then we continue to write $\La$ for ${}_{\chi}\La$.)  We also have a map $\chi^{-1} \colon G \to \mr{Aut}_{\La}(\La)$ defined by $\chi^{-1}(g)(\lambda) = \lambda \cdot \chi(g^{-1})(1)$ for $g \in G$ and
$\lambda \in \Lambda$,
so the resulting object $\La_{\chi}$ is a $\La[G]$-module.
The relationship between ${}_{\chi}\La$ and $\La_{\chi}$ can be expressed by the evaluation-at-$1$
maps
\begin{eqnarray*}
    \Hom_{\La}(\La_{\chi},\La) \lra {}_{\chi}\La &\mr{and}&
    \Hom_{\Lo}({}_{\chi}\La,\La) \lra \La_{\chi},
\end{eqnarray*}
which are isomorphisms of $\Lo[G]$ and $\La[G]$-modules, respectively.

\bl \label{tech chain}
    Suppose that $\La$ is a flat $R$-algebra, and let $\chi \colon G \to \mr{Aut}_{\Lo}(\La)$
    be a homomorphism.
        For any bounded above complex $A$ of $R[G]$-modules that are finitely presented
    over $R$ and
        any bounded below complex $B$ of $R$-modules, the two maps
        \begin{eqnarray*}
                &\theta \colon {}_{\chi}\La \ot_R \Hom_R(A,B) \lra
        \Hom_{\La}(\La_{\chi} \ot_R A,\La\ot_RB)& \\
                &\theta' \colon
        \La_{\chi} \ot_R \Hom_R(A,B) \lra \Hom_{\Lo}({}_{\chi}{\La} \ot_R A, \La \ot_R B)&
        \end{eqnarray*}
    defined by $\theta(\la\ot f)(\mu\ot x) = \mu\la\ot f(x)$ and $\theta'(\la \ot f)(\mu \ot x)
    = \la\mu \ot f(x)$
        are isomorphisms of complexes of $\Lo[G]$-modules and of complexes of $\La[G]$-modules,
    respectively.
\el

\bpf
    We focus on the case of $\theta$, as the other case follows from it.
    Since ${}_{\chi}\La$ is $R$-flat,
    Lemma \ref{tech 2} (with $\Om = \Lo[G]$ and $\Xi = R[G]$)
    implies that the natural map
    $$
        {}_{\chi}\La \ot_R \Hom_R(A,B) \lri \Hom_R(A,{}_{\chi}\La \ot_R B)
    $$
    is an isomorphism of complexes of
    $\Lo[G] \ot_R R[G]$-modules
    if we take the actions of $G$ independently,
    and hence of $\Lo[G]$-modules if we take the
    $G$-actions prescribed earlier in this subsection.
    The inverse of the isomorphism provided by the evaluation map
    provides the first of the
    isomorphisms of complexes of $\Lo[G]$-modules
    $$
        {}_{\chi}\La \ot_R B \lri \Hom_{\La}(\La_{\chi},\La) \ot_R B
        \lri \Hom_{\La}(\La_{\chi}, \La \ot_R B),
    $$
    the canonical second map being an isomorphism as $\La_{\chi}$ is free of rank $1$ over $\La$.
    We therefore have an isomorphism
    $$
        \Hom_R(A,{}_{\chi}\La \ot_R B) \lri \Hom_R(A,\Hom_{\La}(\La_{\chi}, \La \ot_R B)).
    $$
    Finally, the inverse of the adjoint morphism is an isomorphism
    $$
        \Hom_R(A,\Hom_{\La}(\La_{\chi}, \La \ot_R B)) \lri
        \Hom_{\La}(\La_{\chi} \ot_R A, \La \ot_R B)
    $$
    of complexes of $\Lo[G]$-modules,
    and the resulting composite of three isomorphisms is easily computed to be $\theta$.
\epf

\subsection{Projective modules over a profinite ring} \label{projmods}
Let $\La$ be a profinite ring, and fix a directed fundamental system
$\mc{I}$ of open neighborhoods of zero consisting of two-sided
ideals of $\La$. We say that a topological $\La$-module $M$ is
endowed with the $\mc{I}$-adic topology if the collection $\{\A
M\}_{\A\in\mc{I}}$ forms a fundamental system of neighborhoods of
zero.  It was shown in \cite[Section 3.1]{Lim-adic} that any
finitely generated compact (Hausdorff) $\La$-module necessarily has
the $\mc{I}$-adic topology, and, moreover, any homomorphism between
such modules is necessarily continuous.

In this subsection, we recall several facts about projective
$\La$-modules that will be of use to us.  We denote the abelian
category of compact $\La$-modules by $\C_{\La}$. The free profinite
$\La$-module on a set $X$ is canonically isomorphic to the
topological direct product of one copy of $\La$ for each element of
$X$ \cite[Proposition 7.4.1]{Wil}, and a profinite $\La$-module $P$
is projective if and only if it is continuously isomorphic to a
direct summand of the free profinite module on a set of generators
of $P$ \cite[Proposition 7.4.7]{Wil}. In particular, the category
$\C_{\La}$ has enough projectives. Any projective object in
$\C_{\La}$ that is finitely generated over $\La$ is a projective
$\La$-module, and conversely, any finitely generated projective
$\La$-module endowed with the $\mc{I}$-adic topology is a projective
object in $\C_{\La}$ \cite[Proposition 3.1.8]{Lim-adic}.

Recall
that the projective dimension of
an abstract $\La$-module $M$ is the minimum integer $n$ (if it
exists) such that there is a resolution of $M$ by projective
$\La$-modules
\[ 0\lra P^{-n} \lra \cdots \lra P^{-1} \lra P^0 \lra M \lra 0.\]
The topological projective dimension
of a compact $\La$-module is defined similarly, replacing projective
$\La$-modules by projective objects in $\C_{\La}$.
The notions of projective dimension and topological projective dimension
coincide on finitely generated (compact) $\La$-modules.

For the remainder of the subsection, we suppose that $\La$ is left
Noetherian. Note that any projective resolution $P$ of a finitely
generated $\La$-module $M$ is quasi-isomorphic to a projective
resolution $Q$ of $M$ by finitely generated $\La$-modules via a map
$Q \to P$ compatible with the augmentations to $M$. In particular,
the projective dimension of $M$ is the length of its shortest
resolution by projectives in the category of finitely generated
$\La$-modules.

In general, if $\mf{C}$ is a category with objects that are
$\La$-modules, we will use $\mf{C}^{\La-\ft}$ to denote the full
subcategory of objects that are finitely generated (i.e., of finite
type) over $\La$. Since $\Mod^{\La-\ft}_{\La}$ has enough
projectives, the equivalent category $\C^{\La-\ft}_{\La}$ also has
enough projectives.
Denote by $\DD^-_{\La-\ft}(\Mod_{\La})$ the full subcategory of
    the bounded above derived category
$\DD^-(\Mod_{\La})$ which has as its objects those bounded above
complexes $X$ of $\La$-modules for which all of the $H^i(X)$ are
finitely generated $\La$-modules.  The following standard lemma
tells us that any such complex is quasi-isomorphic to a complex of
finitely generated modules. (See \cite[Proposition 3.2.6]{Ne} for an
analogous statement, which has a similar proof.)

\bl \label{complex finite cohom}
    Let $\Om$ be a left Noetherian ring.
    Every bounded above complex $X$ of $\Om$-modules for which every
    $H^i(X)$ is a finitely generated $\Om$-module has a quasi-isomorphic subcomplex
    of finitely generated $\Om$-modules.
\el

In sum, we
have equivalences of categories
\[ \DD^-(\C^{\La-\ft}_{\La}) \lri
\DD^-(\Mod^{\La-\ft}_{\La}) \lri
\DD^-_{\La-\ft}(\Mod_{\La}), \]
the first being induced by the forgetful functor and the second
by the inclusion of categories $\Mod^{\La-\ft}_{\La}\hra \Mod_{\La}$.
We use these equivalences to identify these categories with each other.

\subsection{Completed tensor products} \label{ctp}
In this subsection, we review some basic facts about completed tensor
products and briefly study their derived functors. Let $R$ be a
commutative profinite ring, and let $\La$, $\Om$, and $\Si$ be
profinite $R$-algebras, by which we shall mean that they are
profinite and the maps from $R$ to their centers are continuous.

Let $\mc{I}$ (resp., $\mc{J}$) denote a directed fundamental
system of open neighborhoods of zero consisting of two-sided ideals
of $\La$ (resp., $\Om$). We then define a completed tensor
product algebra by
$$
    \Om \cotimes{R} \Lo = \plim_{\A\in\mc{I},\ \mf{B}\in\mc{J}}\,
    \Om/\mf{B} \ot_R (\La/\mf{A})^{\circ}.
$$
This is clearly a profinite $R$-algebra.

We shall denote the category of compact $\Om \cotimes{R}
\Lo$-modules by $\C_{\Om-\La}$. Let $M$ be an object of
$\C_{\Om-\La}$, and let $N$ be an object of $\C_{\La-\Si}$. By
\cite[Lemma 5.1.1(a)]{RZ}, we have that the set of open, finite index
$\Om \cotimes{R} \Lo$-submodules (resp., $\La \cotimes{R}
\So$-submodules) of $M$ (resp., $N$) forms a basis of
neighborhoods of zero on $M$ (resp., $N$). We define the completed
tensor product to be the compact $\Om \cotimes{R}
\So$-module
\[
    M\cotimes{\La}N = \plim_{M',N'}\,M/M'\ot_{\La}N/N',
\]
where $M'$ and $N'$ run through the respective bases for $M$ and
$N$, with the topology defined by the inverse limit. The completed
tensor product is associative and commutative (i.e., $M
\cotimes{\La} N \cong N \cotimes{\Lo} M$) in the same sense as the
usual tensor product.

Note that the canonical $\Lambda$-balanced map
$$
        t \colon M\times N \lra M\cotimes{\La}N
$$
induces a homomorphism
$$
    M \ot_{\La} N \lra M \cotimes{\La} N
$$
of $\Om \ot_R \So$-modules that has dense image. The
completed tensor product of $M$ and $N$ then satisfies the following
universal property (see \cite[Section 5.5]{RZ} in the case that
$\Om = \Si = R$).

\bl \label{univprop}
    For any compact $\Om \cotimes{\La} \So$-module $L$ and any continuous,
    $\La$-balanced, left $\Om$-linear and right $\Si$-linear
    map $f \colon M\times N \to L$, there is a unique continuous map
    $\hat{f} \colon M\cotimes{\La}N \to L$ of $\Om \cotimes{\La} \So$-modules
    such that $\hat{f} \circ t = f$.
\el

It follows that, in defining the completed tensor
product, it suffices to run through a basis of neighborhoods of zero
consisting of open $\Om \cotimes{R} \Lo$-submodules of $M$ and a basis of
neighborhoods of zero consisting of open $\La \cotimes{R} \So$-submodules of $N$.

The following is also standard.

\bl \label{completed tensor product} Let $M$ and $N$ be objects of
$\C_{\Om-\La}$ and $\C_{\La-\Si}$, respectively.
\begin{enumerate}
    \item[$(a)$] Suppose that $M = \plim M_{\alpha}$ and $N=\plim N_{\beta}$, where each
    $M_{\alpha}$ $($resp., $N_{\beta})$ is a compact $\Om \cotimes{R} \Lo$-module
    $($resp., compact $\La \cotimes{R} \So$-module$)$. Then there is an isomorphism
    \[ M\cotimes{\La}N \cong \plim_{\alpha,\beta} M_{\alpha} \cotimes{\La} N_{\beta}  \]
    of compact $\Om \cotimes{R} \So$-modules.
    \item[$(b)$] The canonical map $M \ot_{\La} N \to M\cotimes{\La}N$
    is an isomorphism if either $M$ is
    finitely generated as a $\Lo$-module or $N$ is finitely generated as a
    $\La$-module.
    \item[$(c)$] The functor
    \[M\cotimes{\La}- \colon \C_{\La-\Si} \lra \C_{\Om-\Si}\]
    is right exact.
\end{enumerate}
\el

The next lemma describes connections between projective
objects in different categories of compact modules and flat
objects in categories of abstract modules.
A ring is said to be left (resp., right) coherent
if all its finitely generated left (resp., right) ideals are finitely presented.

\bl \label{projobj}\ \begin{enumerate}
        \item[$($a$)$] If $\La$ is left $($resp., right$)$ coherent,
        then every projective object in $\C_{\Lo}$ $($resp., $\C_{\La})$
        is flat with respect to the usual left $($resp., right$)$ tensor product over $\La$.
        \item[$($b$)$] If $\Om$ $($resp.,\ $\La)$ is a projective object of $\C_R$,
        then every projective object in $\C_{\Om-\La}$ is a projective object in $\C_{\Lo}$
        $($resp.,\ $\C_{\Om})$.
\end{enumerate}
\el

\begin{proof}
    (a)
    Since every projective object in $\C_{\Lo}$ is continuously isomorphic to a direct summand of a
    direct product of copies of $\Lo$, it suffices to
    show that any direct product of copies of $\La$ is flat as an abstract $\Lo$-module.
    By a theorem of Chase (cf.\ \cite[Theorem 4.47]{Lam}), this is equivalent to the fact that $\La$ is left
    coherent.

    \vspace{2ex}
    (b) Suppose that $\Om$ is projective in $\C_R$, so $\Om$ is continuously isomorphic to
    a direct summand of a
    direct product of copies of $R$.
    Lemma \ref{completed tensor product}(a) then tells us that
    $\Om \cotimes{R} \Lo$ is topologically a direct summand of a direct product of copies of
    $\Lo$, hence is projective.
\end{proof}

Let $G$ be a profinite group. We use $\C_{\La,G}$ to denote the
category of compact $\La$-modules with a continuous commuting action
of $G$, the morphisms being continuous homomorphisms of $\La[G]$-modules.
Note the following.

\begin{remark}
    The category $\C_{\La,G}$ is equivalent to the category $\C_{\La\ps{G}}$,
    where $\La\ps{G}$ is given the profinite topology defined by
    $$
        \La\ps{G} \cong \plim_{\mf{A} \in \mc{I}} \plim_{N} (\La/\mf{A})[G/N],
    $$
    where $N$ runs over the open normal subgroups of $G$.
\end{remark}

Consequently, $\C_{\La,G}$ is an abelian category with enough
projectives, and every element of $\C_{\La,G}$ is an inverse limit of
finite $\La[G]$-quotients (see also \cite[Section
3.2]{Lim-adic}).  Moreover, $\La\ps{G}$ is a projective object of
$\C_{\La}$ by \cite[Lemma 5.3.5(d)]{RZ}, so the forgetful functor
$\C_{\La,G} \to \C_{\La}$ takes projective objects to projective
objects. To shorten notation, we use $\C_{\Om-\La,G}$ to denote the
category $\C_{\Om \cotimes{R} \Lo,G}$.

We remark that if $M$ and $N$ are, respectively, objects of $\C_{\Om-\La,G}$
and $\C_{\La-\Si,G}$, possibly with trivial $G$-actions,  we
may give $M \cotimes{\La} N$ the structure of an object of
$\C_{\Om-\Si,G}$ via the diagonal action of $G$. That is, the
$G$-action is defined by choosing bases of open $(\Om \cotimes{R}
\Lo)[G]$-submodules of $M$ (resp., $(\La \cotimes{R}
\So)[G]$-submodules of $N$) and taking the inverse limits of the
tensor products of the finite quotients.

\begin{remark}
    The analogous result to Lemma \ref{completed tensor product} holds, as a consequence of
    said lemma, if we take
    $M$ and $N$ to be objects of $\C_{\Om-\La,G}$ and $\C_{\La-\Si,G}$, respectively.
    (That is, in part (a), one must take $M_{\alpha}$ and $N_{\beta}$ to be objects of these categories
    to attain an isomorphism in $\C_{\Om-\Si,G}$, and in part (c), the functor is now a
    functor from $\C_{\La-\Si,G}$ to $\C_{\Om-\Si,G}$.)
\end{remark}

Note that we may form completed tensor products of bounded above
complexes, or of a bounded complex with any complex,
as with the usual tensor products.

\bp \label{derived comp tensor}
    If $\Om$ is projective in $\C_R$, then the completed tensor product induces the following
    derived bifunctor
    \[
        - \dotimes{\La} - \colon \DD^-(\C_{\Om-\La,G}) \times \DD^-(\C_{\La-\Si,G})
                \lra \DD^-(\C_{\Om-\Si,G}),
    \]
    where $A\cotimes{\La}B$ represents $A \dotimes{\La} B$ if the terms of $A$
    are projective as objects in $\C_{\Lo}$. Furthermore, there is a commutative diagram
    \[
        \SelectTips{cm}{} \xymatrix@C=1in{
            \DD^-(\C_{\Om-\La,G}) \times \DD^-(\C_{\La-\Si,G}) \ar[d]
            \ar[r]^-{- \dotimes{\La} -} & \DD^-(\C_{\Om-\Si,G}) \ar[d] \\
        \DD^-(\C_{\Lo,G}) \times \DD^-(\C_{\La-\Si,G}) \ar[r]^-{- \dotimes{\La} -}
        & \DD^-(\C_{\So,G})
        }
    \]
    in which the vertical arrows are induced by forgetful functors, and if $\La$ is
    left Noetherian, there is a commutative diagram
    \[
        \SelectTips{cm}{} \xymatrix@C=1in{
        \DD^-(\C_{\Om-\La,G}) \times \DD^-(\C_{\La-\Si,G})
            \ar[r]^-{- \dotimes{\La} -} & \DD^-(\C_{\Om-\Si,G}) \ar[dd] \\
        \DD^-(\C_{\Om-\La,G}) \times \DD^-(\C_{\La-\Si,G}^{\La-\ft}) \ar[u] \ar[d] \\
        \DD^-(\Mod_{(\Om \ot_R \Lo)[G]})  \times \DD^-(\Mod_{(\La \ot_R \So)[G]})
        \ar[r]^-{- \ot_{\La}^{\Li} -} & \DD^-(\Mod_{(\Om \ot_R \So)[G]})
        }
    \]
    in which the vertical arrows are induced by forgetful functors and embeddings of
    categories.
\ep

\bpf
    The proof of the first part and
    the commutativity of the first set of diagrams follow by similar arguments
    to those of Propositions \ref{derived Hom}
    and  \ref{derived tensor G}, making use of Lemma \ref{projobj}(b)
    (using bounded above complexes of projective objects in place
    of q-projective complexes).  The commutativity
    of the final diagram then
    follows from Lemmas \ref{completed tensor product}(b) and  \ref{projobj}.
\epf

Let  $\mr{Ch}^b_{\Lo-\mr{perf}}(\C_{\Om-\La,G})$ denote the category
of bounded complexes in $\C_{\Om-\La,G}$ that are quasi-isomorphic
to bounded complexes in $\C_{\Om-\La,G}$ of objects that are projective in $\C_{\Lo}$,
and let $\DD^b_{\Lo-\mr{perf}}(\C_{\Om-\La,G})$ denote its derived category.

\bp \label{derived comp tensor bdd}
    If $\Om$ is projective in $\C_R$, then we have
    a derived bifunctor
    $$
        - \dotimes{\La} - \colon \DD^b_{\Lo-\mr{perf}}(\C_{\Om-\La,G}) \times \DD(\C_{\La-\Si,G})
                  \lra \DD(\C_{\Om-\Si,G}),
    $$
    where $A \cotimes{\La} B$ represents $A \dotimes{\La} B$ if $A$ is a bounded complex
    with terms that are projectives
    in $\C_{\Lo}$.
 \ep

 \begin{proof}
    Let $B$ be a complex of objects of $\C_{\Om-\La,G}$.
    Since every object of 
    $\mr{Ch}^b_{\Lo-\mr{perf}}(\C_{\Om-\La,G})$ is by definition quasi-isomorphic to
    a bounded complex of objects that are acyclic for the functor $- \cotimes{\La} B$,
    we have a derived functor $- \dotimes{\La} B$.
    Let $A$ be a bounded complex in $\C_{\Om-\La,G}$
    of projectives in $\C_{\Lo}$ and $f \colon B \to B'$ be a quasi-isomorphism.
    Then $f$ induces isomorphisms between the $E_2$-terms of the convergent spectral sequence
    $$
        E_2^{i,j}(B) = H^i(A \cotimes{\La} H^j(B)) \Rightarrow H^{i+j}(A \cotimes{\La} B)
    $$
    and its analogue for $B'$, and therefore it induces isomorphisms on the abutments.
\end{proof}

\subsection{Ind-admissible modules}

The notion of an ind-admissible $R[G]$-module
was introduced in \cite[Section 3.3]{Ne} for a complete commutative Noetherian
local ring $R$ with finite residue field and a profinite group $G$.
An $R[G]$-module is ind-admissible if it can be
written as a union of $R[G]$-submodules that are finitely generated over $R$
and on which $G$ acts continuously with respect to the topology defined by
the maximal ideal of $R$.  In this section, we discuss an analogous construction
of ind-admissible modules over noncommutative profinite rings.

As in Nekov\'a\v{r}'s treatment, we do not consider the seemingly
delicate issue of placing topologies on ind-admissible modules, as
it proves unnecessary.  In particular, it is still possible to
define the continuous cochain complex of an ind-admissible module as
a direct limit. Like the category of discrete modules,
and as opposed to the category of compact modules,
the category of ind-admissible modules has enough injectives.  Thus,
the cochain functor on bounded below complexes of ind-admissible modules
induces a functor on the derived category of such complexes that is computable using injective resolutions.

We maintain the notation of Section \ref{ctp}.  Moreover, we suppose
that $\Om$ is left Noetherian and $\La$ is right Noetherian. For an
$(\Om \ot_R \Lo)[G]$-module $M$, we denote by $\Sa(M)$ the set of
$(\Om \ot_R \Lo)[G]$-submodules of $M$ that are finitely generated
as $\Lo$-modules and on which $\Om$ and $G$ act continuously with
respect to the $\mc{I}$-adic topology. The following is a
straightforward generalization of \cite[Lemma 3.3.2]{Ne}.

\bl \label{ind-ad1} Let $M$ be an $(\Om \ot_R \Lo)[G]$-module.
\begin{enumerate}
    \item[$(a)$] If $M' \in\Sa(M)$, then $N\in\Sa(M)$
    for every $(\Om \ot_R \Lo)[G]$-submodule $N$ of $M'$.
    \item[$(b)$] If $f \colon M\rightarrow N$ is a homomorphism of $(\Om \ot_R \Lo)[G]$-modules
    and $M' \in\Sa(M)$,  then $f(M')\in \Sa(N)$.
    \item[$(c)$] If $M'$, $M'' \in\Sa(M)$, then
    $M'+M''\in\Sa(M)$.
\end{enumerate}
\el

\begin{proof} For part (a), Corollary 3.1.6 and Proposition 3.1.7 of \cite{Lim-adic} imply that the
    subspace topology on $N$ from the $\mc{I}$-adic topology on $M'$
    agrees with the $\mc{I}$-adic topology on $N$, which implies that
    $\Om$ and $G$ act continuously
    on $N$.  In (b), the continuity of the
    $\Om$ and $G$-actions
    on the finitely generated $\Lo$-module
    $f(M')$ is a consequence of the fact that
    the map $M' \to f(M')$ is a continuous quotient map with respect to
    the $\mc{I}$-adic topology (as follows from Corollary 3.1.5 and
    Proposition 3.1.7 of \cite{Lim-adic}).  Part
    (c) follows from (b), using the addition map $M \times M \to M$.
\end{proof}

Note that Lemma \ref{ind-ad1}(c) implies that $\Sa(M)$ is a directed
set with respect to inclusion.
We say that an $(\Om \cotimes{R} \Lo)[G]$-module $M$ is (right) ind-admissible if
$$
    M = \bigcup_{N \in \Sa(M)} N.
$$
We list some basic properties of the full subcategory $\II_{\Om-\La,G}$ of
$\Mod_{(\Om \ot_R \Lo)[G]}$ with objects the ind-admissible $(\Om \ot_R \Lo)[G]$-modules.

\bl \label{ind-ad main}\
    \begin{enumerate}
        \item[$(a)$] The category $\II_{\Om-\La,G}$ is abelian and stable under
        subobjects, quotients and colimits.
        \item[$(b)$] The embedding functor
        $$
            i \colon \II_{\Om-\La,G} \hookrightarrow \Mod_{(\Om \ot_R \Lo)[G]}
        $$
        is exact and
        is left adjoint to the functor
        $$
            j \colon \Mod_{(\Om \ot_R \Lo)[G]} \ra \II_{\Om-\La,G}
        $$
        that takes a module $M$ to the union of the elements of $\Sa(M)$.
        \item[$(c)$]
        The category
        $\II_{\Om-\La,G}$ has enough injectives.
        \item[$(d)$] Let $M$ be an ind-admissible $(\Om \ot_R \Lo)[G]$-module, and let $N$ be a
        finitely generated $\Lo$-submodule
        of $M$.
        Then $(\Om \ot_R \Lo)[G]\cdot N$ is an
        ind-admissible $(\Om \ot_R \Lo)[G]$-module which is a finitely generated $\Lo$-module.
        \item[$(e)$]  Let $M$ be an $(\Om \ot_R \Lo)[G]$-module. Then $M\in\Sa(M)$ if and only
        if $M$ is an ind-admissible $(\Om \ot_R \Lo)[G]$-module which is finitely generated
        as a $\Lo$-module.
    \end{enumerate}
\el

\bpf
    For parts (a), (b) and (c), similar arguments to those of \cite[Proposition
    3.3.5]{Ne} apply.   The ``only if" direction of (e) is obvious.  For (d),
    since $N$ is $\Lo$-finitely generated, we can find a finite subset
    $\{ M_1, \ldots ,M_n \}$ of
    $\Sa(M)$ such that
    \[  N \sbs M_1+\cdots + M_n.\]
    The assertion then follows from Lemma \ref{ind-ad1}(c) and  the ``only if" direction of (e).
    The ``if" direction of (e) follows from (d), since $M
    = (\Om \ot_R \Lo)[G]\cdot M$.
\epf

If $G$ is trivial, we write $\II_{\Om-\La}$ for $\II_{\Om-\La,G}$.  We leave to the reader the proof of the following.

\bl \label{ind-ad ind} \
    \begin{enumerate}
        \item[$(a)$] Let $M$ be an $(\Om \ot_R \Lo)[G]$-module.  If $N$ is an
        $(\Om \ot_R \Lo)[G]$-submodule of $M$, then $N \in \Sa(M)$ if and only if
        $N$ is an object of $\C_{\Om-\La,G}^{\Lo-\ft}$.
        \item[$(b)$] The category $\II_{\Om-\La,G}$ is equivalent to
        the ind-category of $\C_{\Om-\La,G}^{\Lo-\ft}$.
        \item[$(c)$]
        The category $\II_{\Om-\La}$
        is equivalent to the full subcategory of $\II_{\Om-\La,G}$ with objects
        the modules on which $G$ acts trivially.
    \end{enumerate}
\el

\begin{remark}
    In the case that $\Om = R$, we denote $\II_{R-\La,G}$ by $\II_{\Lo,G}$, and the
    latter category is equivalent to the ind-category of $\C_{\Lo,G}^{\Lo-\ft}$.
    The subcategory $\II_{\Lo}$ is simply $\Mod_{\Lo}$.
\end{remark}

\bl \label{cts homs ft}
    Let $A$ be an object of $\C_{\Om,G}^{\Om-\ft}$, and let $B$ be an object of
    $\C_{\Om-\La,G}^{\Lo-\ft}$.
    Then $\Hom_{\Om}(A,B)$ with the $\mc{I}$-adic topology
    is an object of $\C_{\Lo,G}^{\Lo-\ft}$.
\el

\begin{proof}
    Since $A$ is finitely generated over $\Om$ and $B$ is a compact $\Om$-module, we have by
    \cite[Lemma 3.1.4(3)]{Lim-adic} that $\Hom_{\Omega}(A,B) = \Hom_{\Omega,\cts}(A,B)$,
    where the latter group is the group of continuous homomorphisms
    of $\Omega$-modules.
    Note that $\Lo$
    acts continuously on $\Hom_{\Om}(A,B)$ with respect to the compact-open topology
    by \cite[Proposition 3(a)]{F}, and similarly $G$ acts continuously on it
    as a consequence of \cite[Lemma 2]{F}.
    We have a continuous isomorphism of $\Lo[G]$-modules
    $$
        \Hom_{\Om}(A,B) \lri \plim_{\beta} \Hom_{\Om}(A,B_{\beta}),
    $$
    where $B_{\beta}$ runs over the finite Hausdorff
    $(\Om \cotimes{R} \Lo)[G]$-quotients of $B$.
    As $B_{\beta}$ is finite and therefore discrete, the group $\Hom_{\Om}(A,B_{\beta})$
    is finite and discrete as well, so $\Hom_{\Om}(A,B)$ is compact.

    Note that $\Hom_{\Om}(A,B)$ injects into $\Hom_{\Om}(\Om^r,B) \cong B^r$ for some $r \ge 0$,
    since $A$ is $\Om$-finitely generated.  We therefore have that $\Hom_{\Om}(A,B)$ is
    $\Lo$-finitely generated.   Finally,
    \cite[Proposition 3.1.7]{Lim-adic} then implies that the compact-open topology
    on $\Hom_{\Om}(A,B)$ agrees with the $\mc{I}$-adic topology.
\end{proof}

We define
$$
    \Hom_{\Om,\cts}(-,-) \colon \C_{\Om,G} \times \II_{\Om-\La,G} \lra \II_{\Lo,G}  .
$$
by
$$
    \Hom_{\Om,\cts}(A,B) = \ilim_{\al}
    \ilim_{\beta} \Hom_{\Om}(A_{\alpha},B_{\beta}),
$$
with $A = \plim A_{\alpha}$, where $A_{\alpha}$ runs over the
quotients of $A$ in
$\C_{\Om,G}^{\Om-\ft}$, and $B = \ilim B_{\beta}$, where $B_{\beta}$
runs over the elements of $\Sa(B)$, i.e., by Lemma \ref{ind-ad ind},
the subobjects of $B$ that lie in $\C_{\Om-\La,G}^{\Lo-\ft}$.
Note that it actually suffices to let $A_{\alpha}$ and $B_{\beta}$ run over cofinal subsets.

\begin{remark} \label{limit cts hom}
    We note that
    $$
        \Hom_{\Om,\cts}(A,B) \cong \ilim_{\beta}  \Hom_{\Om,\cts}(A,B_{\beta}),
    $$
    where $\Hom_{\Om,\cts}(A,B_{\beta})$ here is the $(\Om \ot_R \Lo)[G]$-module of continuous
    $\Om$-module homomorphisms in the usual sense.
\end{remark}

We also note the following.

\bl \label{compare Hom}
    Let $A$ be an object of $\C_{\Om,G}^{\Om-\ft}$, and let $B$ be an object of
    $\II_{\Om-\La,G}$.  Then we have
    $$
        \Hom_{\Om,\cts}(A,B) = \Hom_{\Om}(A,B).
    $$
\el

\begin{proof}
    This follows from the computation
    $$
        \Hom_{\Om,\cts}(A,B) \cong  \ilim_{\beta} \Hom_{\Om,\cts}(A,B_{\beta})
        = \ilim_{\beta} \Hom_{\Om}(A,B_{\beta}) \cong \Hom_{\Om}(A,B),
    $$
    the latter isomorphism following from the case in which $A$ is free of finite rank, since $A$
    is finitely presented over $\Om$.
\end{proof}

As usual, we may extend our definition to consider the complex of continuous
homomorphisms from a complex of compact modules to a complex of ind-admissible modules,
supposing that at least one of these complexes is bounded above and
at least one is bounded below.

\bp \label{derived cts Hom}
    There is a derived bifunctor
    $$
        \RHom_{\Om,\cts}(-,-) \colon \DD^-(\C_{\Om,G})^{\circ} \times
        \DD^+(\II_{\Om-\La,G}) \lra \DD^+(\II_{\Lo,G}),
    $$
    and $\RHom_{\Om,\cts}(A,B)$ can be represented by $\Hom_{\Om,\cts}(A,B)$ if
    the terms of $A$ are projective as objects of $\C_{\Om}$.
    Moreover, we have a commutative diagram
    $$
        \SelectTips{cm}{} \xymatrix@C=1in{
        \DD^-(\C_{\Om,G})^{\circ} \times \DD^+(\II_{\Om-\La,G}) \ar[r]^-{\RHom_{\Om,\cts}(-,-)}
        & \DD^+(\II_{\Lo,G}) \ar[dd] \\
        \DD^-(\C_{\Om,G}^{\Om-\ft})^{\circ} \times \DD^+(\II_{\Om-\La,G})
        \ar[u] \ar[d] \\
        \DD^-(\Mod_{\Om[G]})^{\circ} \times \DD^+(\Mod_{(\Om\ot_R\Lo)[G]}) \ar[r]^-{\RHom_{\Om}(-,-)}
        & \DD^+(\Mod_{\Lo[G]})   }
    $$
    in which the vertical arrows arise from forgetful functors and
    the natural embeddings of abelian categories.
\ep

\begin{proof}
    The existence of the derived bifunctor in question
    follows from an argument similar to that of Proposition \ref{derived Hom},
    and the second statement is a consequence of the fact
    that a projective object in $\C_{\Om,G}$ is a projective object of
    $\C_{\Om}$, noting Remark \ref{limit cts hom}.

    Let $A$ be a bounded above
    complex of objects of $\C_{\Om,G}^{\Om-\ft}$, and let $B$ be a bounded below complex
    of objects of $\II_{\Om-\La,G}$.  Choose bounded above complexes $L$, $P$ and
    $Q$, consisting of projective objects in $\C_{\Om,G}$, $\Mod_{\Om[G]}$
    and $\Mod_{\Om}^{\Om-\ft}$ (i.e., $\C_{\Om}^{\Om-\ft}$)
    respectively, with quasi-isomorphisms to $A$ of complexes
    in the respective categories. Projectivity yields natural maps
    $Q \to L$, $Q \to P$ and $P \to L$
    of complexes in $\C_{\Om}$, $\Mod_{\Om}$ and $\Mod_{\Om[G]}$,
    respectively. We then have a commutative diagram
    $$
        \SelectTips{cm}{} \xymatrix{
        \Hom_{\Om,\cts}(L,B) \ar[r] \ar[d] & \Hom_{\Om,\cts}(Q,B) \ar[d] \\
        \Hom_{\Om}(P,B) \ar[r] & \Hom_{\Om}(Q,B).
        }
    $$
    Here, the right-hand
    vertical arrow is the identity map by Lemma \ref{compare Hom}
    (viewing $Q$ as a complex of objects in $\C_{\Om,G}^{\Om-\ft}$ with trivial $G$-action),
    and the lower horizontal map is a quasi-isomorphism as a consequence of Proposition
    \ref{derived Hom G}, as $\Om[G]$ is $\Om$-free.
    We want to show that the left-hand
    vertical arrow is a quasi-isomorphism of complexes of $\Lo[G]$-modules, and we
    will be done if we can show the upper horizontal arrow is a quasi-isomorphism of
    complexes of $\Lo$-modules.

    In the case that $B$ is a module, exactness of direct limit reduces us to the
    case that $B$ is an object of $\C_{\Om-\La,G}^{\Lo-\ft}$, hence of $\C_{\Om}$.
    In this case, since both $L$ and $Q$ are complexes of projectives
    in $\C_{\Om}$, we obtain that the map is a quasi-isomorphism.

    In the general setting, since $L$ and $Q$ are bounded above and $B$ is bounded below,
    we have
    $$
        \Hom_{\Om,\cts}^n(X,B) = \bigoplus_{j \in \Z} \Hom_{\Om,\cts} (X^j, B^{j+n})
    $$
    for $X = L$ and $X = Q$, and the direct sum commutes with direct
    limits.  As $B$ is the direct limit of its truncations
    $\tau_{\leq i}B$ (with $i$th terms $\ker d_B^i$),
    we may therefore assume that $B$ is bounded with $B^j = 0$ for $j > i$.
    In this case, brutal truncations (i.e., $\sigma_{\le i} B$ having $i$th term $B^i$)
    provide an exact triangle
    \[  (\sigma_{\leq i-1}B)[i-1] \lra B^i \lra (\sigma_{\leq i}B)[i] \lra (\sigma_{\le i-1} B)[i], \]
    and we are reduced recursively to the above-proven case of a module.
\end{proof}

The following analogue of Lemma \ref{tech chain} will be of later use.

\bl \label{tech chain indad}
    Suppose that $\La$ is a flat, Noetherian $R$-algebra, and let
            $\chi \colon G \to \mr{Aut}_{\Lo}(\La)$ be a homomorphism.
        Let $A$ be a bounded above complex of objects in $\C_{R,G}$, and let $B$
        be a bounded below complex of $R$-modules.  Then there are isomorphisms
        \begin{eqnarray*}
                    &\theta \colon {}_{\chi} \La \ot_R \Hom_{R,\cts}(A,B) \lra
                \Hom_{\La,\cts}(\La_{\chi} \cotimes{R} A,\La\ot_RB)& \\
                    &\theta' \colon
                \La_{\chi} \ot_R \Hom_{R,\cts}(A,B) \lra
                \Hom_{\Lo,\cts}({}_{\chi}{\La} \cotimes{R} A, \La \ot_R B)&
        \end{eqnarray*}
    that arise as direct limits of the maps in Lemma \ref{tech chain}.
\el

\begin{proof}
    Write $A = \plim A_{\alpha}$ with $A_{\alpha}$ an object of  $\C_{R,G}^{R-\ft}$ and
    $B = \ilim B_{\beta}$ with $B_{\beta}$ an object of $\C_R^{R-\ft}$.  Note that
    $$
        \La_{\chi} \cotimes{R} A \cong \plim_{\al}\, (\La_{\chi} \ot_R A_{\alpha}),
    $$
    so, applying Lemma \ref{tech chain}, we have
    \begin{align*}
        \Hom_{\La,\cts}(\La_{\chi} \cotimes{R} A, \La \ot_R B) &=
        \ilim_{\al,\beta} \Hom_{\La}(\La_{\chi} \ot_R A_{\al},\La \ot_R B_{\beta}) \\
        &\cong \ilim_{\al,\beta} {}_{\chi}\La \ot_R \Hom_R(A_{\al},B_{\beta})\\
        &\cong {}_{\chi} \La \ot_R \Hom_{R,\cts}(A,B).
    \end{align*}
    This yields the first isomorphism, and the argument for the second is similar.
\end{proof}

Let $M$ be an ind-admissible
$(\Om \ot_R \Lo)[G]$-module. The continuous cochain complex of
$G$ with values in $M$ is the complex
of $\Om \ot_R \Lo$-modules that is the direct limit of
complexes
\[
    C(G,M) = \ilim_{N \in\Sa(M)}C(G,N),
\]
where $C(G,N)$ is the usual complex of (inhomogeneous) continuous cochains,
with $N$ given the $\mc{I}$-adic topology.

We remark that, in this definition,
it clearly suffices to take the direct limit over a cofinal subset of $\Sa(M)$.  In particular,
if $M$ itself is finitely generated over $\Lo$, then by Lemma \ref{ind-ad main}(e) the above definition of
$C(G,M)$ as a direct limit agrees with its definition as continuous cochains, considering
$M$ as an object of $\C_{\Om-\La,G}$.  Finally, we can extend the above definition
to consider the total cochain complex of a complex $M$
of ind-admissible modules, which has $k$th term
$$
    C^k(G,M) = \bigoplus_{i+j = k} C^i(G,M^j)
$$
and differentials as in \cite[(3.4.1.3)]{Ne}, and the cochain functor induces a functor
$$
    \R\Ga(G,-) \colon \DD^+(\II_{\Om-\La,G}) \lra \DD^+(\Mod_{\Om \ot_R \Lo})
$$
between derived categories.  Since $\II_{\Om-\La,G}$ has enough injectives,
it can always be computed on a bounded below complex using a
quasi-isomorphic complex of injectives.

\section{Cohomology groups of induced modules}

In this paper, we are interested in the case that our profinite ring
$\La$ is the Iwasawa algebra of a compact $p$-adic Lie group.  With
this in mind, we introduce some notation that will be used from this
point forward.  Fix a prime $p$. Let $R$ be a commutative complete
Noetherian local ring with maximal ideal $\m$ and residue field $k$,
where $k$ is finite of characteristic $p$.  We let $\Gamma$ denote a
compact $p$-adic Lie group. We are interested in the completed group ring
$\La = R\ps{\Ga}$.

We note that $\La$ is a profinite ring, endowed
with the topology given by the canonical isomorphism
$$
    \La \lri \plim_{U \in \U} \plim_{n \ge 1} R/\m^n[\Gamma/U],
$$
where $\U$ denotes the set of open normal subgroups of $\Gamma$.
In fact, it is a projective object in $\C_R$ (cf.\ \cite[Lemma 5.3.5(d)]{RZ}).
Moreover, we have the following (cf. \cite[Theorem 8.7.8]{Wil}).

\setcounter{thm}{0}

\bp \label{Noetherian compact ring} The ring $\La$ is Noetherian. \ep

\bpf
    As $R$ is a complete local Noetherian ring with finite residue field $k$, the Cohen
    structure theorem \cite[Theorem 12]{Coh} implies that it
    is isomorphic to a quotient of a power series ring $S$ in $n = \dim_k \m/\m^2$
    variables over the ring of Witt vectors $\mc{O}$ of $k$.
    Therefore, $R\ps{\Ga}$ is a quotient
    of $S\ps{\Ga}$, so it suffices to prove that $S\ps{\Ga}$ is Noetherian.  Since
    $S\ps{\Ga} \cong \mc{O}\ps{\Zp^n\times \Ga}$ and $\Zp^n \times \Ga$ is a compact $p$-adic
    Lie group, we have that $S\ps{\Ga}$ is Noetherian by a mild extension of a classical
    theorem of Lazard's (cf.\ \cite[Corollary 2.4]{Ven05} and
    \cite[Proposition V.2.2.4]{Laz}).
\epf

Since Noetherian rings are necessarily coherent, we may
apply Lemma \ref{projobj}(a) to conclude from Proposition \ref{Noetherian compact ring}
that $\La$ is a flat $R$-algebra.

\subsection{Induced modules and descent} \label{indmod}

In this subsection, we will exhibit an interesting spectral sequence relating
the cohomology of an induced module over an Iwasawa algebra to the
cohomology of the module itself.  We start work in a more general setting.

Let $G$ be a profinite group.  Let $\Si$ be a left coherent and right Noetherian
profinite $R$-algebra.  For a complex $M$ of objects in $\C_{\Si,G}$, we let $C(G,M)$ denote its
complex of continuous $G$-cochains, we let $\R\Ga(G,M)$ denote the corresponding
object in $\DD(\Mod_{\Si})$, and we let $H^i(G,M)$ denote its $i$th continuous 
$G$-cohomology (or, more precisely, hypercohomology) group.
As in \cite[Proposition 3.2.11]{Lim-adic}, the functor
$$
    \R\Ga(G,-) \colon \DD^+(\C_{\Si,G}) \lra \DD^+(\Mod_{\Si})
$$
is well-defined and exact.
(While $\C_{\Si,G}$ may not have enough injectives, the treatment of \cite{KS}
asserts the existence of $\DD^+(\C_{\Si,G})$ and $\DD(\C_{\Si,G})$
under certain set-theoretic assumptions, which we make here.)

In the case that $G$ has finite $p$-cohomological dimension, we can do better.
For a double complex $X$ (or, by abuse of notation, its total complex) and $n \in \Z$,
we let $\tau_{\ge n}^{II}(X)$ (resp., $\tau_{\le n}^{II}(X)$)
be the total complex of the quotient complex (resp., subcomplex)
of $X$ with $j$th row equal to $\tau_{\ge n}(X^{\bullet j})$
(resp., $\tau_{\le n}(X^{\bullet j})$).  
For lack of a sufficiently precise reference, we provide a short proof of
the following. 

\bl \label{hypercoh spec seq}
    Suppose that $G$ has finite $p$-cohomological dimension.  Then we have a convergent
    hypercohomology spectral sequence
    $$
        E_2^{i,j} = H^i(G,H^j(M)) \Rightarrow H^{i+j}(G,M)
    $$
    for any complex $M$ of objects in $\C_{\Si,G}$.
\el

\begin{proof}
	Since $G$ has finite cohomological dimension, for sufficiently large $k$
	we have a convergent hypercohomology spectral sequence
	$$
		E_2^{i,j} = H^i(G,H^j(M)) \Rightarrow H^{i+j}(\tau_{\le k}^{II} C(G,M))
	$$
	arising from the filtration on rows of the indicated truncation of the 
	double $G$-cochain complex of $M$ (cf.\ \cite[\S II.2]{NSW}).
	As $C(G,M)$ is the direct limit of the complexes $\tau_{\le k}^{lI} C(G,M)$
	under maps that are quasi-isomorphisms by the above convergence,
	the natural map 
	$$
		H^{i+j}(\tau_{\le k}^{II} C(G,M)) \to H^{i+j}(G,M)
	$$
	is an isomorphism by exactness of the direct limit.
\end{proof}

As quasi-isomorphisms of complexes induce isomorphisms on the $E_2$-terms of their hypercohomology spectral sequences, we have the following corollary.

\bc \label{derived functor fin cd}
    Suppose that $G$ has finite $p$-cohomological dimension.
    The functor $C(G,-)$ preserves quasi-isomorphisms of chain complexes
    and induces an exact functor
    $$
        \R\Ga(G,-) \colon \DD(\C_{\Si,G}) \to \DD(\Mod_{\Si}).
    $$
\ec

Fix a profinite $R$-algebra $\Om$ that is a projective object of $\C_R$
for the remainder of the section.
We now derive the following spectral sequence.

\bp \label{tensor cohomology p-adic group}
        	Let $Y$ be a complex of objects in $\C_{\Si,G}$, and let $N$ be a bounded above
        	complex of objects in $\C_{\Om-\Si}^{\So-\ft}$.  Consider the conditions
        	\begin{enumerate}
                	\item[$\mr{(1)}$] $G$ has finite $p$-cohomological dimension,
        		\item[$\mr{(2)}$] $N$ is bounded with terms of finite projective dimension over $\So$.
    	\end{enumerate}
    	If $\mr{(1)}$ holds and $Y$ is bounded above, $\mr{(2)}$ holds and $Y$ is bounded below, or both
    	$\mr{(1)}$ and $\mr{(2)}$ hold, then we have an isomorphism
            $$
                    N\ot_{\Si}^{\Li}\R\Ga(G, Y) \lri \R\Ga(G, N \dotimes{\Si} Y)
            $$
            in $\DD(\Mod_{\Om})$.
\ep

\bpf
        Let $L$ be a bounded above complex of projective objects in $\C_{\Om-\Si}$ mapping
        quasi-isomorphically to $N$.  If (1) holds and $Y$ is bounded above, we set $P = L$, and we have
    by Proposition \ref{derived comp tensor} that $P \cotimes{\Si} Y$ represents
    $N \dotimes{\Si} Y$ in $\DD^-(\C_{\Om,G})$.
    Otherwise, there exists a quasi-isomorphic bounded quotient complex
    $P = \tau_{\ge n} L$ of $L$ of objects in
    $\C_{\Om-\Si}$, through which the quasi-isomorphism $L \to N$ factors, such that
    the terms of $P$ are projective in $\C_{\So}$.
    Proposition \ref{derived comp tensor bdd}
    then tells us that $P \cotimes{\Si}Y$
    represents $N \dotimes{\Si} Y$ in $\DD(\C_{\Om,G})$.
    In all cases, the terms of $P$ are flat $\So$-modules by Lemma \ref{projobj},
    so Proposition \ref{derived tensor} tells us that $P\ot_{\Si}C(G, Y)$ represents
    $N\ot_{\Si}^{\Li}\R\Ga(G,Y)$ in $\DD(\Mod_{\Om})$.

    We have a map of complexes of $\Om$-modules
    \[
        P\ot_{\Si}C(G, Y) \lra
        C(G, P \cotimes{\Si} Y)
    \]
    with sign conventions as in \cite[Proposition 3.4.4]{Ne}.
    (The continuity of the cochains in the image is insured, for instance, by the fact that
    any term of $P$ is a topological direct summand of a direct product of copies of $\Si$.)
    It suffices to show that this map is a quasi-isomorphism.

    Let $Q$ be a bounded above complex of finitely generated projective
    $\So$-modules mapping quasi-isomorphically to $N$, which we take to be bounded if (2) holds.
    We then have a map of
    complexes of $R$-modules
    \[
        Q \ot_{\Si}C(G, Y) \lri
        C(G, Q \ot_{\Si} Y)
    \]
    that, much as in \cite[Proposition 3.4.4]{Ne}, is an isomorphism.
    (To see that it is an isomorphism, note that it is immediate if $Q$ is a finitely generated free
    $\So$-module, and therefore also for direct summands of such, i.e., the finitely generated
    projective $\So$-modules.  The case of a complex is then immediate.)
    By the projectivity of $Q$ in both the categories of abstract and compact
    $\So$-modules, we have a commutative diagram
    $$
        \SelectTips{cm}{} \xymatrix{
            Q \ot_{\Si}C(G, Y) \ar[r]^{\sim} \ar[d] & C(G, Q \ot_{\Si} Y) \ar[d] \\
            P\ot_{\Si}C(G, Y) \ar[r] & C(G, P \cotimes{\Si} Y)
        }
    $$
    in which the vertical arrows are quasi-isomorphisms of complexes of $R$-modules,
    noting Lemma \ref{completed tensor product}(b) and, when (1) holds,
    Corollary \ref{derived functor fin cd}.
    It follows that the lower horizontal arrow is a quasi-isomorphism as well.
\epf

We make a couple of remarks.

\begin{remark} \label{tensor coh fin gp}
        Suppose that $G$ is a finite group.  The proof of Proposition \ref{tensor cohomology p-adic group}
    under condition (2) then goes through word-for-word for bounded $Y$,
        with Tate cochains and the resulting derived complexes replacing
        the usual cochains of $G$ and its derived complexes.
\end{remark}

\begin{remark}
    In \cite[Proposition 1.6.5(3)]{FK}, an identical isomorphism to that of
    Proposition \ref{tensor cohomology p-adic group} is proved for
    a very general class of ``adic" rings $\Om$ and $\La$, with stronger conditions on $G$,
    $N$, and $Y$.
\end{remark}

From now on, we focus on case that $\Si$ is the Iwasawa algebra $\La = R\ps{\Ga}$.
Fix a continuous homomorphism
$\chi \colon G \to \Ga$ of profinite groups.  Since $\Ga$ may be
viewed as a subgroup of $\mr{Aut}_{\Lo}(\La)$ by left
multiplication, we may define the $\Lo[G]$-module ${}_{\chi}\La$ as
in Section \ref{mod gp rings}. If $A$ is a
$\La \ot_R \Lo$-module, we then set ${}_{\chi}A = {}_{\chi}\La \ot_{\La} A$
and $A_{\chi} = A \ot_{\La} \La_{\chi}$, which are $\Lo[G]$ and
$\La[G]$-modules, respectively.

Let $M$ be a $R[G]$-module.  As in \cite[Section 5.1]{Lim-adic}, we
define a $\La[G]$-module $\F_{\Ga}(M)$ by

\[
    \F_{\Ga}(M) = \plim_{U\in\U} (R[\Gamma/U]_{\chi}\ot_R M)
\]
with $G$ acting diagonally and $\La$ acting on the left on the terms
in the inverse limit. This construction applied to compact modules
provides a functor
$$
    \F_{\Ga} \colon \C_{R,G} \lra \C_{\La,G},
$$
which is exact as $R[\Gamma/U]$ is $R$-flat and the inverse limit is
exact on inverse systems in $\C_{\La,G}$. We extend this to an exact functor on complexes in the
obvious fashion.

Given an complex $T$ of objects in $\C_{R,G}$, we have, noting Lemma \ref{completed tensor product}(a), a natural continuous isomorphism
$$
    \F_{\Ga}(T) \cong \La_{\chi} \cotimes{R} T
$$
of complexes in $\C_{\La,G}$.  We may therefore use Lemma
\ref{completed tensor product}(b) to note that $\F_{\Ga}$ takes a
complex $T$ of objects in in $\C_{R,G}^{R-\ft}$ to the complex
$\La_{\chi} \ot_R T$ of objects in $\C_{\La,G}^{\La-\ft}$. We make
the latter identification freely.

We now supply a key ingredient for descent.

\bl \label{bbc}
        	Let $T$ be a complex in $\C_{R,G}$.
        	Let $\Gamma'$ be a quotient of $\Gamma$ by a closed normal subgroup, and set
   	$\La' = R\ps{\Gamma'}$.  Let $N$ be a bounded above complex of objects of $\C_{\Om-\La'}$.
    	Then $N$ can be viewed as a complex of objects in $\C_{\Om-\La}$
        	via the quotient map $\pi \colon \La \tha \La'$.
    	Suppose either that $T$ is bounded above or that $N$ is bounded and its
    	terms have finite projective dimension in $\C_{\Lo}$ and in $\C_{(\La')^{\circ}}$.
    	Then $\pi$ induces an isomorphism
        	\[
                	N \dotimes{\La} \F_{\Ga}(T) \lri N \dotimes{\La'} \F_{\Gamma'}(T)
        	\]
        	in $\DD^-(\C_{\Om,G})$.
\el

\bpf
        Let $Q \to N$ be a quasi-isomorphism, where $Q$ is a bounded above complex of
        $\C_{\Om-\La'}$-projectives, and note that its terms are projective objects of
        $\C_{(\La')^{\circ}}$ by Lemma \ref{projobj}(b).
        If the terms of $N$ have finite projective dimension in $\C_{(\La')^{\circ}}$, then for sufficiently
        small $n$, the truncation $\tau_{\ge n} Q$  also consists of projectives in $\C_{(\La')^{\circ}}$,
        so we may assume in that case that $Q$ is bounded.  Next, and similarly, let $P \to Q$ be a
        quasi-isomorphism in $\C_{\Om-\La}$ with $P$ a bounded above complex in
        $\C_{\Om-\La}$ of projectives in $\C_{\Lo}$, and suppose that $P$ is bounded
        if the terms of $N$ have finite projective dimension in $\C_{\Lo}$.
        The derived completed tensor products in question are then
            represented by $P \cotimes{\La} \La_{\chi} \cotimes{R} T$
        and $Q \cotimes{\La'} \La'_{\chi} \cotimes{R} T$, respectively, by Propositions
        \ref{derived comp tensor} and \ref{derived comp tensor bdd}.  Since the induced map
        $P \cotimes{\La} \La_{\chi} \to Q \cotimes{\La'} \La'_{\chi}$ is clearly a quasi-isomorphism
        and the terms of these complexes are projective in $\C_R$, we have by the same
        propositions the desired quasi-isomorphism between the two complexes.
\epf

We are particularly interested in the case that the complex $N$ above is the Iwasawa algebra of a 
quotient of $\Gamma$, for which lemma the following is crucial.

\bl \label{fin proj dim}
    	Let $\Ga_0$ be a closed normal subgroup of $\Ga$ with no elements of order $p$.
    	Set $\Ga' = \Ga/\Ga_0$ and $\Lambda' = R\ps{\Ga'}$.
        	Then $\La'$ has finite projective dimension over $\La$.
\el

\begin{proof}
    	Let $\cTor_i^{\Xi}(-,-)$ denote the $i$th derived bifunctor of the completed tensor product
    	over a profinite $R$-algebra $\Xi$.  A standard argument yields a
    	convergent spectral sequence
        	$$
            	H_i(\Ga, \cTor_j^{R}(\La',Z)) \Rightarrow \cTor_{i+j}^{\La}(\La',Z)
        	$$
        	for any compact $\La$-module $Z$, as in \cite[\S V.2, Exercise 3]{NSW}.
        	Since $\La'$ is a projective object in $\C_R$, the spectral sequence degenerates
        	to yield isomorphisms
        	$$
            	H_n(\Ga, \La'\cotimes{R}Z) \cong \cTor_{n}^{\La}(\La',Z).
       	$$
        	Moreover, by a standard extension of Shapiro's lemma 
	\cite[Theorem 6.10.9]{RZ}, we have
    	$$
        		H_n(\Ga_0,Z) \cong H_n(\Ga,\La' \cotimes{R} Z).
    	$$
    	Since $\Ga_0$ has finite $p$-cohomological dimension by
        	\cite[Corollaire 1]{Se}, there then exists $n_0 \ge 0$ independent of $Z$
    	such that $\cTor_{n}^{\La}(\La',Z) = 0$ for every $n > n_0$.
        	It follows from \cite[Proposition 5.2.11]{NSW}
        	that $\La'$ has finite topological projective dimension as a compact $\La$-module, so it has
        	finite projective dimension over $\La$ by the discussion of Section \ref{projmods}.
\end{proof}

The following descent spectral sequence is reduced by the above results to a special case of
Proposition \ref{tensor cohomology p-adic group}.

\bt \label{descseq}
    Let $T$ be a complex in $\C_{R,G}$.
    Let $\Ga'$ be a quotient of $\Gamma$ by a closed normal subgroup $\Ga_0$,
    and set $\Lambda' = R\ps{\Gamma'}$.  Consider the conditions
        \begin{enumerate}
            	\item[$(1)$] $G$ has finite cohomological dimension, 
		\item[$(2)$] $\Ga_0$ contains no elements of order $p$.
        \end{enumerate}
        Suppose that $(1)$ holds and $T$ is bounded above, $(2)$ holds and $T$ is
        bounded below, or both $(1)$ and $(2)$ hold.
    Then we have an isomorphism
    \[
        \La' \ot_{\La}^{\Li}\R\Ga(G, \F_{\Ga}(T)) \lri \R\Ga(G, \F_{\Ga'}(T))
    \]
    in $\DD(\Mod_{\La'})$.
\et

\begin{proof}
    This is simply Proposition \ref{tensor cohomology p-adic group} for
    $N = \La'$ and $Y = \F_{\Ga}(T)$, noting Lemma \ref{fin proj dim} and
    applying the isomorphism of Lemma \ref{bbc}.
\end{proof}

\subsection{Finite generation of cohomology groups} \label{fingen}

In this subsection, we shall show that the cohomology groups of
induced modules are finitely generated under a certain assumption on
the group.  We maintain the notation of the previous subsection.

\bl \label{fg La R}
    Suppose that $\Ga$ is pro-$p$.  Let $M$ be a compact $\La$-module.
    Then $M$ is finitely generated over $\La$ if and only if $R\ot_{\La}M$
    is finitely generated over $R$.
\el

\bpf
    Since $\Gamma$ is pro-$p$, the Jacobson radical $\M$ of $\Lambda$
    is $\m\La + I$, where $I$ is the augmentation ideal of $\Lambda$
    (see \cite[Proposition 5.2.16(iii)]{NSW}).  This implies that
    $$
        M/\M M \cong M_{\Gamma}/\m M_{\Gamma},
    $$
    where $M_{\Gamma} = M/I M$.
    Therefore, Nakayama's lemma tells us that $M_{\Gamma}$ is finitely generated
    over $R$ if and only if $M/\M M$ is finite.
    On the other hand, Nakayama's lemma for compact $\La$-modules (cf.\ \cite[Lemma 5.2.18(ii)]{NSW})
    tells us that $M$ is finitely generated over $\La$ if and only if $M/\M M$ is finite.
\epf

Lemma \ref{fg La R} requires a compact $\La$-module.
Therefore, we give a sufficient condition
for an abstract $\La$-module to be a compact $\La$-module
under an appropriate topology.

\bl \label{fg La R 2} Suppose that $M$ is an abstract $\La$-module
which is the inverse limit of an inverse system of finite quotient modules.
Then $M$ is a compact
$\La$-module with respect to the
resulting profinite topology.
\el

\bpf We need only to show that the $\La$-action
\[ \theta \colon \La\times M \lra M\]
is continuous with respect to the topology given by the inverse
limit. Suppose first that $M$ is finite. Let $\M$ denote the
Jacobson radical of $\La$. Then $\M^n M$ stabilizes, and it follows
from Nakayama's lemma that we have $\M^n M =0$ for large
enough $n$.  By \cite[Corollary 5.2.19]{NSW},
there exist $r \ge 0$ and $U\in\U$ such that
$\m^r\La + I(U)\sbs \M^n$, where $I(U)$ denotes the ideal of $\La$ generated by
the augmentation ideal in $R\ps{U}$.
Hence, $\theta$ is continuous with respect to the discrete topology on $M$.

In general, we can write $M \cong \plim M/M_{\al}$, where $\{ M_{\al} \}$
is a directed system of $\La$-submodules of finite index.
Let $(\la, x)\in \theta^{-1}(y + M_{\al})$ for
$\la\in \La$ and $x,y\in M$.
Since $M/M_{\al}$ is finite,
it follows from the above discussion that there exist $r$ and
$U$  such that
$$
    (\la + \m^r\La + I(U))\cdot (x+M_{\al})\sbs y + M_{\al},
$$
as desired. \epf

We shall also require the following lemma.

\bl \label{fg tor} If $M$ is a finitely generated $\La$-module, then
$\Tor_{i}^{\La}(R,M)$ is finitely generated over $R$ for every $i \ge 0$.  \el

\begin{proof}
    To see this, we first
    choose a resolution $P$ of $M$ consisting of finitely generated projective
    $\La$-modules. Then $R \ot_{\La}P$ is a complex of finitely generated
    $R$-modules.  Therefore, the homology groups $\Tor_{i}^{\La}(R,M)$ of the latter complex
    are finitely generated over $R$.
\end{proof}

We recall from \cite[Proposition 4.2.3]{Ne} that if $G$ is
a profinite group such that $H^{i}(G,M)$ is finite for every finite
$G$-module $M$ of $p$-power order
and $i \ge 0$, then $H^{i}(G,T)$ is a finitely
generated $R$-module for every
$T\in\C_{R,G}^{R-\ft}$ and $i\geq 0$.
We will refer to a profinite group $G$ as $p$-cohomologically finite if there exists
an integer $n$ such that, for every finite $G$-module $M$ of $p$-power order,
$H^i(G,M) = 0$ for all $i > n$ and
$H^i(G,M)$ is finite for all $i$.
We remark that if $G$ is $p$-cohomologically finite, then so is every open subgroup of $G$, as follows directly from Shapiro's lemma.

\bt \label{fg cohomology}
    Suppose that $G$ is $p$-cohomologically finite.
    Let $T$ be a complex of objects in $\C_{R,G}^{R-\ft}$.
    Then the cohomology groups $H^i(G, \F_{\Ga}(T))$ are finitely generated over
    $\La$ for all $i$.  That is, $\R\Ga(G,\F_{\Ga}(T))$ is an object of $\DD_{\La-\ft}(\Mod_{\La})$.
\et

\bpf
    Suppose for now that $T$ is concentrated in degree $0$.
        Let $\Gamma_0$ be an open uniform normal pro-$p$ subgroup of
        $\Gamma$, and set $H = \chi^{-1}(\Ga_0)$.
        Fix a set of double coset representatives $\gamma_1, \ldots, \gamma_t$ of
        $\Ga_0\backslash \Ga/ \chi(G)$.
        For each $i$, define $\chi_i \colon G \to \Ga$ by $\chi_i(g) = \gamma_i g \gamma_i^{-1}$.
     The reader may check that
    \begin{eqnarray*}
            &\displaystyle \bigoplus_{i=1}^t \Big(\Z[G] \ot_{\Z[H]}
            R\ps{\Ga_0}_{\chi_i}\Big)\stackrel{\sim}{\lra} \La_{\chi}& \\
            &( g\ot \la_0)_i \mapsto \la_0 \ga_i \chi(g)^{-1}&
    \end{eqnarray*}
        is an isomorphism of $R\ps{\Ga_0}[G]$-modules.  Given this, Shapiro's lemma
    induces isomorphisms
        $$
            H^j(G, \F_{\Ga}(T)) \cong \displaystyle \bigoplus_{i=1}^t
            H^j(H, R\ps{\Ga_0}_{\chi_i}\ot_R T)
        $$
        of $R\ps{\Ga_0}$-modules.
        The finite
        generation of $H^j(G, \F_{\Ga}(T))$ over $\La$ then reduces to
        the finite generation of each $H^j(H,
        R\ps{\Ga_0}_{\chi_i}\ot_R T)$ over $R\ps{\Ga_0}$.
        Thus, we can and do assume that $\Gamma$ has no elements of finite order.

    Since $R$ has finite projective dimension over $\La$ by Lemma \ref{fin proj dim},
    Theorem \ref{descseq} provides the convergent spectral sequence
    \[
        E_2^{r,s} = \Tor_{-r}^{\La}(R,H^s(G,\F_{\Ga}(T))) \Rightarrow H^{r+s}(G, T).
    \]
    Since $G$ has finite $p$-cohomological dimension, say $n$,
    we may suppose inductively that $H^i(G, \F_{\Ga}(T))$ is a finitely generated
    $\La$-module for all $i$ greater than some $j \le n$.  Note that $E_{\infty}^{0,j}$
    is a quotient of $E_2^{0,j}$ and that
    $H^j(G,T)$ is a finitely generated $R$-module by
    \cite[Proposition 4.2.3]{Ne}.  Since
    $E_{\infty}^{0,j}$ is a subquotient of $H^{j}(G,T)$,
    it is also a finitely generated $R$-module.  On the other hand, it follows from the
    definition of $E_{\infty}^{0,j}$ that the kernel of the surjective map
    $E_2^{0,j} \to E_{\infty}^{0,j}$ is isomorphic to a subquotient of the finite direct sum
    \[
        \bigoplus_{i = j+1}^n E_2^{j-i-1,i}
        = \bigoplus_{i = j+1}^n \mr{Tor}^{\La}_{i-j+1}(R, H^i(G,\F_{\Ga}(T))).
    \]
    By our induction hypothesis and Lemma \ref{fg tor}, the above
    module is a finitely generated $R$-module.
    It follows that
    $$
        E_2^{0,j} \cong R \otimes_{\Lambda} H^j(G,\F_{\Ga}(T))
    $$
    is a finitely generated $R$-module.
    As $H^j(G,\F_{\Ga}(T))$ is an inverse limit of finite $\La$-modules by
    \cite[Proposition 5.2.4]{Lim-adic}, Lemma \ref{fg La R 2}
    gives it the structure of a compact $\La$-module, and we may therefore apply
    Lemma \ref{fg La R} to conclude that $H^j(G,\F_{\Ga}(T))$ is finitely generated over $\La$.

        Now let us work without condition on $\Gamma$ or the complex $T$.
    By Lemma \ref{hypercoh spec seq} and the exactness of $\F_{\Ga}(-)$, we have the
    convergent spectral sequence of $\La$-modules
    \[
        H^r(G,\F_{\Ga}(H^s(T))) \Rightarrow H^{r+s}(G,\F_{\Ga}(T)).
        \]
    Thus, $H^i(G,\F_{\Ga}(T))$ has a filtration consisting of subquotients of the
    finitely generated $\La$-modules $H^j(G ,\F_{\Ga}(H^{i-j}(T)))$ with $0 \le j \le n$.
        Hence, $H^i(G,\F_{\Ga}(T))$ is also a finitely generated $\La$-module.
\epf

\begin{remark}
        In the case that $\Gamma$ is abelian, Theorem \ref{fg cohomology}
    for bounded below $T$ is essentially \cite[Proposition 4.2.3]{Ne},
        as $R\ps{\Gamma_0}$ for an open pro-$p$, torsion-free subgroup $\Gamma_0$ is
        itself a commutative complete local Noetherian ring with finite residue field.
        In the case that $T$ is a bounded complex of $R$-projectives, it is a corollary of
        \cite[Proposition 1.6.5(2)]{FK}.
\end{remark}

\subsection{Duals of induced modules} \label{duals}

In this subsection, we describe the pairings and the resulting cup
products that will be used in proving a duality theorem for the
cohomology of induced modules.  Key to our discussion is the
following lemma, which we record for convenience.  We let $\Om$ and
$\Si$ denote auxiliary Iwasawa algebras: e.g., $\Om = R\ps{\Phi}$
for some compact $p$-adic Lie group $\Phi$.

\bl \label{totcupprod}
    Suppose that $M$ is a bounded above
    complex of objects of $\C_{\Om-\La,G}$, let $N$ be
    a bounded above complex of objects of $\C_{\La-\Si,G}$, and let $L$ be a
    bounded above complex of objects of $\C_{\Om-\Si,G}$.  Then any map
    $$
            \phi \colon M \cotimes{\La} N \lra L
    $$
    of complexes in $\C_{\Om-\Si,G}$ gives rise to a cup product morphism
    $$
            \cup \colon C(G,M) \ot_{\La} C(G,N) \lra C(G,L)
    $$
    of complexes of $\Om \ot_{\La} \So$-modules.
\el

\begin{proof}
    Briefly, the point is that to give such a map $\phi$ is equivalent to producing a
    collection of continuous,
    $\La$-balanced, $G$-equivariant,
    left $\Om$-linear, and right $\Si$-linear
    pairings
    $$
        \langle -,- \rangle_{mn} \colon M^m \times N^n \lra L^{m+n}
    $$
    that are compatible with coboundaries in the sense of \cite[Section 3.3]{Lim-adic}.
    The cup product $\cup$ then arises from the induced cup products
    $$
        \cup_{mn}^{ij} \colon C^i(G,M^m) \ot_{\La} C^j(G,N^n) \lra C^{i+j}(G,L^{m+n})
    $$
    by combining them to a map of the total complexes with appropriate signs
    as in \cite[(3.4.5.2)]{Ne}: $\cup = ((-1)^{in}\cup_{mn}^{ij})$.
\end{proof}

Let us fix some notation.  Let $\iota \colon \La \ra \La$ denote the
unique continuous $R$-algebra homomophism that satisfies
$\iota(\gamma) = \gamma^{-1}$ for all $\gamma \in \Gamma$.  If $A$
is a $\La[G]$-module, we let $A^{\iota}$ denote the $\Lo[G]$-module
that is $A$ as an $R[G]$-module but on which any
$\lambda \in \La$ now acts by left multiplication by $\iota(\lambda)$. We extend
this to complexes in the obvious fashion.  Of particular interest to
us is $\F_{\Ga}(T)^{\iota}$ for a complex $T$ of objects in
$\C_{R,G}$.  Since $\iota$ induces a continuous isomorphism
of $\Lo[G]$-modules ${}_{\chi}\Lambda \xrightarrow{\sim} \Lambda_{\chi}^{\iota}$,
we can and will make the identification
$$
    \F_{\Ga}(T)^{\iota} = {}_{\chi} \Lambda \cotimes{R} T,
$$
of complexes in $\C_{\Lo,G}$ from this point forward.

\bl \label{fg pairing1}
        For objects $M$ and $N$ in $\C_{R,G}$,
        the pairing
    \begin{eqnarray*}
        &\langle - , - \rangle \colon \F_{\Ga}(M) \times \F_{\Ga}(N)^{\iota} \lra
            \La\cotimes{R}(M\cotimes{R}N)&\\
            &(\la\ot m , \mu\ot n) \mapsto \la\mu\ot m\ot n&
    \end{eqnarray*}
        is continuous, $R$-balanced,
        $G$-equivariant with respect to the trivial action on $\La$ in the target,
        and is left and right $\La$-linear.
            Moreover, if $M$ and $N$ are bounded above
        complexes of objects in $\C_{R,G}$, we have a morphism
        $$
                \phi \colon \F_{\Ga}(M) \cotimes{R} \F_{\Ga}(N)^{\iota}\lra \La\cotimes{R}(M\cotimes{R}N)
        $$
        of complexes in $\C_{\La-\La,G}$.
\el

\bpf
    We begin by remarking that by Lemma \ref{univprop}, it makes sense to define the pairing
    on pairs of tensors, identified with their images in the completed tensor product.
    The pairing is clearly $R$-balanced, and the reader can check
    the statements on the $\La$, $\Lo$, and $G$-actions.
    The continuity of the pairing reduces immediately to the continuity of multiplication on $\La$.
    The second statement follows easily.
\epf

We next apply these pairings to study a duality between induced modules.
As described in the introduction, the
dualizing complex for $R$ is an object $\w_R$ of $\DD^b_{R-\ft}(\Mod_{R})$
with the property that for every object $M$ of $\DD(\Mod_R^{R-\ft})$, the object
$\RHom_{R}(M,\w_R)$ lies in $\DD_{R-\ft}(\Mod_R)$ and the canonical morphism
\[
    M\lra \RHom_R(\RHom_{R}(M,\w_R), \w_R)
\]
is an isomorphism in $\DD(\Mod_R)$. A dualizing complex exists for
$R$ and is unique up to translation and isomorphism in
$\DD^b_{R-\ft}(\Mod_{R})$ (see \cite[Ch.\ V]{RD}). We choose a
bounded complex $J_R$ of injective $R$-modules which represents a
choice of the dualizing complex in $\DD(\Mod_R)$. (See \cite[Section
0.4]{Ne}, or Section \ref{change of rings} below, for the construction of such a complex.)

Let $T$ be a bounded complex of objects in $\C_{R,G}^{R-\ft}$.
Then $\Hom_{R}(T,J_R)$ is a bounded complex of  ``admissible"
$R[G]$-modules with cohomology groups that are
finitely generated over $R$ (see \cite[(4.3.2)]{Ne}). By \cite[Proposition 3.3.9]{Ne},
there is a subcomplex $T^*$ of
$\Hom_{R}(T,J_R)$ which is a complex of
objects in $\C_{R,G}^{R-\ft}$ (giving said objects the $\m$-adic topology)
and is quasi-isomorphic to $\Hom_{R}(T,J_R)$ via the
inclusion map.

We have a composite morphism
\[
    \pi \colon T\ot_{R}T^{*} \lra T\ot_{R}\Hom_{R}(T,J_R) \lra J_R
 \]
of complexes of $R[G]$-modules, the first morphism being
induced by the inclusion and the second being the usual evaluation map.
By Lemma \ref{fg pairing1}, this in turn induces a composite map
\[
    \overline{\pi} \colon \F_{\Ga}(T) \cotimes{R} \F_{\Ga}(T^*)^{\iota} \stackrel{\phi}{\lra} \La\ot_R T\ot_R
    T^{*} \stackrel{\id\ot\pi}{\lra} \La \ot_R J_R
\]
of complexes of $(\La \ot_R \Lo)[G]$-modules. Twisting
$\overline{\pi}$ by the identity map on an auxilliary complex of
objects of $\C_{\Om-\Lambda}^{\Lo-\ft}$,  Lemma \ref{totcupprod} now
provides the cup product morphisms of the following lemma, which
will be used in the next section.
For a continuous character $\kappa \colon G \to R^{\times}$ and an
$R[G]$-module $M$, we let $M(\kappa)$ denote the $R$-module $M$
with the new commuting $G$-action given by the twist of the original by $\kappa$.
(In addition, $M(\kappa)$ will be taken to maintain any
topology and actions of profinite $R$-algebras with which $M$ may be  endowed.)

\bl \label{cupprod2}
    Let $A$ be a bounded above complex of objects in $\C_{\Om-\La}^{\Lo-\ft}$,
    and let $\kappa \colon G \to R^{\times}$ be a continuous homomorphism.
    The map $\overline{\pi}$ defined above induces a map
    $$
            \overline{\pi}_A \colon A \ot_{\La} \F_{\Gamma}(T) \cotimes{R}
            \F_{\Gamma}(T^*)^{\iota} \lra A \ot_R J_R
    $$
    and, in turn, a well-defined cup product
    $$
        C(G, A \ot_{\La} \F_{\Ga}(T)) \ot_R C(G,\F_{\Ga}(T^*)^{\iota}(\kappa))
        \lra C(G, A \ot_R J_R(\kappa)),
    $$
    which is a map of complexes of $\Om \ot_R \Lo$-modules.
\el

Note that, by definition, the adjoint maps (defined as in
\cite[Lemma 2.2]{Lim-adic})
\begin{eqnarray*}
    \mr{adj}(\pi) \colon T^* \lra \Hom_R(T,J_R) & \text{and} &
    \mr{adj}'(\pi) \colon T \lra \Hom_R(T^*,J_R)
\end{eqnarray*}
of $\pi$ are quasi-isomorphisms, the first being simply the inclusion of complexes
already defined.
Since the terms of $J_R$ are injective $R$-modules, it follows from Proposition
\ref{derived Hom G} that the derived adjoint maps
\begin{eqnarray*}
    \mathbf{adj}(\pi) \colon T^* \lra \RHom_R(T,\omega_R) &\text{and}&
    \mathbf{adj}'(\pi) \colon T \lra \RHom_R(T^*,\omega_R)
\end{eqnarray*}
are isomorphisms in $\DD^b_{R-\ft}(\Mod_{R[G]})$, hence in
$\DD^b_{R-\ft}(\II_{R,G})$, considering the restriction of $\RHom_R(-,-)$ to
a bifunctor
$$
    \DD^-(\C_{R,G}^{R-\ft})^{\circ} \times \DD^+(\Mod_R) \lra \DD^+(\II_{R,G}),
$$
(see \cite[(3.5.9)]{Ne}).

Though we shall not use it later, we feel it important to note that
the derived adjoint maps of $\overline{\pi}$ are also isomorphisms.
While this is easy enough to prove in the bounded below
derived category $\DD^+(\Mod_{\Lo[G]})$
of abstract modules (in the case of $\mathbf{adj}(\overline{\pi})$),
we are interested in $G$-cohomology
groups, so we want such an isomorphism in $\DD^+(\II_{\Lo,G})$.

\bt \label{D(FM)2}
    For any bounded complexes $T$ and $T^*$ of objects in $\C_{R,G}^{R-\ft}$ such that $T^*$ sits
    quasi-isomorphically as a subcomplex of the dual $\Hom_R(T,J_R)$, the derived adjoint maps
    \begin{eqnarray*} &\mathbf{adj}(\overline{\pi}) \colon \F_{\Ga}(T^*)^{\iota} \lra
    \RHom_{\La,\cts}(\F_{\Ga}(T), \La\ot_{R}^{\Li}\w_R)&\\
    &\big(\text{resp., }\mathbf{adj}'(\overline{\pi}) \colon \F_{\Ga}(T) \lra
    \RHom_{\Lo,\cts}(\F_{\Ga}(T^*)^{\iota}, \La\ot_{R}^{\Li}\w_R)\big)&
    \end{eqnarray*}
    are isomorphisms in $\DD^+(\II_{\Lo,G})$ $($resp., $\DD^+(\II_{\La,G}))$.
    Moreover, we have that the derived object
    $\RHom_{\La,\cts}(\F_{\Ga}(T), \La \ot_R^{\Li} \w_R)$
    $($resp., $\RHom_{\Lo,\cts}(\F_{\Ga}(T^*)^{\iota}, \La \ot_R^{\Li} \w_R))$
    can be represented
    by $\Hom_{\La}(\F_{\Ga}(T), \La \ot_R J_R)$ $($resp., $\Hom_{\Lo}(\F_{\Ga}(T^*)^{\iota},\La \ot_R
    J_R))$.
\et

\bpf
    Let $X  = \RHom_{\La,\cts}(\F_{\Ga}(T),\La \ot_R^{\Li} \omega_R)$.
    Since $J_R$ is a complex of $R$-injectives
    and ${}_{\chi} \La$ is flat over $R$,
    Lemma \ref{tech chain indad} implies that the functor $F$ given by
    $$
        \Hom_{\La,\cts}(\La_{\chi} \cotimes{R} -, \La \ot_R J_R)
        \colon \mr{Ch}^-(\C_{R,G})^{\circ} \lra  \mr{Ch}^+(\II_{\Lo,G})
    $$
    is exact.
    In particular, if $P \to T$ is a quasi-isomorphism with $P$ a bounded above
    complex of projective objects
    in $\C_{R,G}$, then $F(T) \to F(P)$ is a quasi-isomorphism.
        As $\La_{\chi} \cotimes{R} P$
    is a complex of objects in $\C_{\La,G}$ that are projective in $\C_{\La}$, we then have that $F(T)$
    represents $X$ by Proposition \ref{derived cts Hom}.

    By Lemmas \ref{completed tensor product}(b) and \ref{compare Hom}, the functor $F$
    agrees with $\Hom_{\La}(\La_{\chi} \ot_R -, \La \ot_R J_R)$
    on $\mr{Ch}^-(\C_{R,G}^{R-\ft})^{\circ}$.
    Therefore, $X$ is represented by
    $\Hom_{\La}(\F_{\Ga}(T), \La \ot_R J_R)$ in $\DD^+(\II_{\Lo,G})$.
    Finally, since the map
    $$T^* \lra \Hom_R(T,J_R)$$
    is a quasi-isomorphism, the adjoint map of complexes
    $$
        {}_{\chi}\La \ot_R T^* \lra \Hom_{\La}(\La_{\chi} \ot_R T, \La \ot_R J_R)
        $$
        is a quasi-isomorphism of complexes of $\Lo[G]$-modules by Lemma \ref{tech chain}.
        Therefore, we have the result for $\mathbf{adj}(\overline{\pi})$, and the proof
    for $\mathbf{adj}'(\overline{\pi})$ is analogous.
\epf

\section{Duality over $p$-adic Lie extensions}

We now turn to arithmetic.  Here, we fix the notation that we shall
use throughout this section. To start, let $p$ be a prime.  We let
$F$ be a global field of characteristic not equal to $p$.  Let $S$
be a finite set of primes of $F$ that, in the case that $F$ is a
number field, contains all primes above $p$ and all real places.
Let $S_f$ (resp., $S_{\RR}$) denote the set of finite places (resp., real places) in $S$.
Let $G_{F,S}$
denote the Galois group of the maximal unramified outside $S$
extension of $F$.  For a place $v$ of $F$, let $F_v$ denote the completion
of $F$ at $v$, and let $G_v$ denote a fixed decomposition group for $v$
in the absolute Galois group of $F$.

We fix a $p$-adic Lie extension $F_{\infty}$ of $F$ that is
unramified outside $S$, and we let $\Gamma$ denote its Galois group.
We let $R$ denote a complete commutative local Noetherian ring with
finite residue field of characteristic $p$, and we set $\La =
R\ps{\Gamma}$.  We take the homomorphism 
$\chi \colon G_{F,S} \to \Ga$ of Section \ref{indmod} to be restriction. 

\subsection{Iwasawa cohomology} \label{fgic}

We recall the following facts, all of which can all be found in
\cite[Chapters VII-VIII]{NSW}. The $G_{F,S}$-cohomology groups of a
finite $G_{F,S}$-module of $p$-power order are all finite. If $p$ is
odd or $F$ has no real places, then $G_{F,S}$ has $p$-cohomological
dimension at most $2$ and so is $p$-cohomologically finite in the
sense of Section \ref{fingen}.  If $v$ is a nonarchimedean place of $F$,
then $G_v$ has $p$-cohomological dimension equal to $2$, and the $G_v$-cohomology
groups of a finite $G_v$-module of $p$-power order are finite as
well.  Of course, the $G_v$-cohomology groups of a finite module are
also finite for archimedean $v$, since $G_v$ is of order dividing
$2$ for such places.

Recall (e.g., from \cite[(5.7.2)]{Ne}) that, for any profinite ring
$\Om$, the (Tate) compactly supported $G_{F,S}$-cochain complex
of a complex of objects $M$ in $\C_{\Om,G_{F,S}}$ is defined as
$$
    C_{(c)}(G_{F,S},M) = \Cone\Bigg( C(G_{F,S},M) \lra \bigoplus_{v \in S_f}
    C(G_v,M)
    \oplus \bigoplus_{v \in S_{\RR}} \widehat{C}(G_v,M)
     \Bigg)[-1],
$$
where $\widehat{C}(G_v,M)$ is defined as in \cite[Section 3.4]{Lim-adic}
and denotes the standard complete complex of Tate $G_v$-cochains for
$M$.

We remark that, in Nekov\'a\v{r}'s notation, $C_{(c)}(G_{F,S},M)$ is
denoted $\widehat{C}_c(G_{F,S},M)$. We use the notation of \cite{FK}, where
parentheses are used to distinguish the latter group from the
compactly supported cochains $C_c(G_{F,S},M)$, for which one uses
the usual cohomology groups at real places.
For archimedean $v$, we will abuse notation and use $\R\Ga(G_v,M)$ to denote the
derived object of the Tate cochains $\widehat{C}(G_v,M)$.
These cochains may as well be taken to be zero for complex places
or for real places if $p$ is odd, $\R\Ga(G_v,M)$ being a zero object in the derived category.
We denote the derived object corresponding to
$C_{(c)}(G_{F,S},M)$ by $\R\Ga_{(c)}(G_{F,S},M)$ and its $i$th
cohomology group by $H^i_{(c)}(G_{F,S},M)$.  By definition, we have
an exact triangle\footnote{We write an exact triangle
    $A \to B \to C \to A[1]$ more compactly as $A \to B \to C$ throughout.}
$$
     \R\Ga_{(c)}(G_{F,S},M) \lra \R\Ga(G_{F,S},M) \lra \bigoplus_{v \in S} \R\Ga(G_v,M)
$$
in $\DD(\Mod_{\Om})$.
The results of and methods used in Section \ref{indmod} allow us to prove
the following descent result for such an exact triangle with induced coefficients.

\bp \label{spectri}
    	Let $\Gamma' = \Gal(F'_{\infty}/F)$ be a quotient of $\Gamma$ by a closed normal
        	subgroup, and set $\La' = R\ps{\Gamma'}$.  Let $T$ be a bounded above
    	complex of objects in $\C_{R,G_{F,S}}$.  Suppose that $p$ is odd, that $F'_{\infty}$ has no 
	real places, or that $T$ is bounded and $F'_{\infty}$ has no real places that become 
	complex in $F_{\infty}$. Then we have an isomorphism of exact triangles
        	$$
          	\SelectTips{cm}{} \xymatrix{
          	\La' \ot_{\La}^{\Li} \R\Ga_{(c)}(G_{F,S},\F_{\Ga}(T)) \ar[r]^-{\sim} \ar[d]
                	& \R\Ga_{(c)}(G_{F,S},\F_{\Ga'}(T)) \ar[d] \\
                	\La' \ot_{\La}^{\Li} \R\Ga(G_{F,S},\F_{\Ga}(T)) \ar[r]^-{\sim} \ar[d] &
        		\R\Ga(G_{F,S}, \F_{\Ga'}(T)) \ar[d] \\
                	\La' \ot_{\La}^{\Li} \bigoplus_{v \in S} \R\Ga(G_v,\F_{\Ga}(T)) \ar[r]^-{\sim}  &
                	\bigoplus_{v \in S} \R\Ga(G_v,\F_{\Ga'}(T))  }
        	$$
    	in $\DD(\Mod_{\La'})$.
\ep

\begin{proof}
    We need only show that the lower and upper horizontal morphisms are isomorphisms.
    The upper morphism exists and is an isomorphism for $T$ bounded above without
    additional assumption.  That is,
    recall that the compactly supported cohomology groups of a module vanish in dimension
    greater than $3$ \cite[Lemma 5.7.3]{Ne}.
    The hypercohomology spectral sequence for compactly supported cohomology therefore
    converges for bounded above complexes, yielding the existence of a functor
    $$
        \R\Ga_{(c)}(G_{F,S},-) \colon \DD^-(\C_{\La,G_{F,S}}) \to \DD^-(\Mod_{\La}).
    $$
    The analogous
    result to Proposition \ref{tensor cohomology p-adic group} then holds by the original argument
    (in the case that (1) holds and $Y$ is bounded above).  That the upper morphism
    exists and is an isomorphism follows immediately as in Theorem \ref{descseq}.

    The summands of the lower morphism corresponding to nonarchimedean $v$ are
    isomorphisms, again for $T$ bounded above without further assumption,
    by Theorem \ref{descseq}.
    For a real place $v$ of $F$, the functor $\R\Ga(G_v,-)$ is well-defined only on
    the bounded derived category.  However, if $v$ becomes complex in $F'_{\infty}$,
    then the relevant hypercohomology spectral sequence implies that the composite
    functor $\R\Ga(G_v,\F_{\Ga'}(-))$ is both well-defined and zero on $\DD(\C_{R,G_v})$.
    So, the $v$-summand of the lower horizontal morphism is trivially an isomorphism.
    If $v$ splits completely in $F_{\infty}/F$, then the modules
    $\La_{\chi}$ and ${}_{\chi} \Lambda$ have trivial $G_v$-actions, so we have isomorphisms
    $\F_{\Ga}(T) \cong \La \ot_R T$ and $\F_{\Ga}(T^*)^{\iota} \cong \La \ot_R T^*$
    of complexes in $\C_{\La,G_v}$ and $\C_{\Lo,G_v}$, respectively.
    Since $G_v$ is finite, the terms of $\widehat{C}(G_v,M)$ for a
    $G_v$-module $M$ are each naturally isomorphic to a finite direct sum of copies of $M$.
    With these identifications, the canonical map
    \begin{eqnarray*}
        \La \ot_R \widehat{C}(G_v,T) \lri \widehat{C}(G_v, \La \ot_R T)
    \end{eqnarray*}
    agrees with the identity map, so is an isomorphism.  Of course, we
    have the corresponding result for $\La'$, and the isomorphism
    $$
        \La' \ot_{\La}^{\Li} (\La \ot_R \widehat{C}(G_v,T)) \lri \La' \ot_R^{\Li}
        \widehat{C}(G_v,T)
    $$
    in $\DD(\Mod_{\La'})$ provides the $v$-summand of the lower isomorphism.
\end{proof}

\begin{remark}
    	Proposition \ref{spectri} holds with ``bounded above" removed if we suppose that
	the kernel of $\Gamma \to \Gamma'$
    	has no elements of order $p$.  This follows quickly if $p$ is odd
    	or $F$ has no real places, but it requires some work
    	in the remaining case that $F'_{\infty}$ has no real places but $F$ does.
    	We omit this for purposes of brevity.
\end{remark}

The next proposition shows in particular that if $T$ is a bounded complex of objects 
of $\C_{R,G_{F,S}}^{R-\ft}$, then the global, local, and
compactly supported cohomology groups of $\F_{\Ga}(T)$ are finitely
generated $\La$-modules.  
(Here, we implicitly identify cochain complexes with their quasi-isomorphic
truncations in the derived category.)

\bp \label{fg type cohomology}
    	Let $T$ be a complex of objects of $\C_{R,G_{F,S}}^{R-\ft}$.
        	\begin{enumerate}
            	\item[$(a)$] The complex $\R\Ga_{(c)}(G_{F,S},\F_{\Ga}(T))$
        		lies in $\DD_{\La-\ft}(\Mod_{\La})$ if $F_{\infty}$ has no real places or $p$ is odd and
        		in $\DD_{\La-\ft}^-(\Mod_{\La})$ if $T$ is bounded above.
        		\item[$(b)$]  The complex $\R\Ga(G_{F,S},\F_{\Ga}(T))$ for $v \in S_f$ lies in
	        $\DD_{\La-\ft}(\Mod_{\La})$ if $F_{\infty}$ has no real places or $p$ is odd and
	        in $\DD_{\La-\ft}^+(\Mod_{\La})$ if $T$ is bounded below.
	        \item[$(c)$] If $T$ is bounded and either $F_{\infty}$ has no real places or $p$ is odd,
	        then $\R\Ga_{(c)}(G_{F,S},\F_{\Ga}(T))$ and $\R\Ga(G_{F,S},\F_{\Ga}(T))$
	        lie in $\DD_{\La-\ft}^b(\Mod_{\La})$.
	        \item[$(d)$] For $v \in S_f$, the complex $\R\Ga(G_v,\F_{\Ga}(T))$
	        lies in $\DD_{\La-\ft}(\Mod_{\La})$ and in $\DD_{\La-\ft}^b(\Mod_{\La})$
	        if $T$ is bounded.
	        \item[$(e)$]  The complex $\R\Ga(G_v,\F_{\Ga}(T))$ for $v \in S_{\RR}$ lies in
	        $\DD_{\La-\ft}(\Mod_{\La})$ if $T$ is bounded.
    \end{enumerate}
\ep

\bpf
    For real $v$, any $i \in \Z$ and $T$ concentrated in degree $0$,
    the Tate cohomology group $\widehat{H}^i(G_v,\F_{\Ga}(T))$ is a $2$-torsion
    subquotient of $\F_{\Ga}(T)$, as $G_v$ has order $2$.  Hence, it is finitely generated
    over $\La$ as $\F_{\Ga}(T)$ is.  The hypercohomology spectral sequence
    allows one to pass to the case of a bounded
    complex $T$, and therefore $\R\Ga(G_v,\F_{\Ga}(T))$ lies in $\DD_{\La-\ft}(\Mod_{\La})$,
    hence (e).  Recall also that $\R\Ga(G_v,\F_{\Ga}(T))$ is a zero object for all complexes $T$
    if $v$ extends to a complex place of $F_{\infty}$.

    If $v \in S$ is nonarchimedean, then $G_v$ is $p$-cohomologically finite by our above remarks, and
    hence Theorem \ref{fg cohomology} implies that $\R\Ga(G_v,\F_{\Ga}(T))$ lies in
    $\DD_{\La-\ft}(\Mod_{\La})$.  In the case that $p$ is odd or $F_{\infty}$ has no real places, Theorem
    \ref{fg cohomology} again implies that $\R\Ga(G_{F,S},\F_{\Ga}(T))$ lies in
    $\DD_{\La-\ft}(\Mod_{\La})$.  (For this, note that if $p = 2$ and $F$ has a real place that
    becomes complex in $F_{\infty}$, it will already have become complex in any extension
    with Galois group $\Gamma/\Gamma_0$, where $\Gamma_0$ is an open uniform pro-$p$
    subgroup of $\Gamma$.)
    In this case, it follows from the above-mentioned exact triangle that
    $\R\Ga_{(c)}(G_{F,S},\F_{\Ga}(T))$ is in $\DD_{\La-\ft}(\Mod_{\La})$.  If $T$ is bounded,
    then clearly all of the above complexes will additionally lie in the bounded derived category.
    In particular, we have (c) and (d) and part of each of (a) and (b).

    The analogue of Theorem \ref{fg cohomology} for compactly
    supported cohomology holds for bounded above $T$ since
    $H^i_{(c)}(G_{F,S},\F_{\Ga}(T))$ vanishes
    for sufficiently large $i$.  The proof is essentially identical but uses
    the semilocal version of Shapiro's lemma \cite[(8.5.3.2)]{Ne}
    in the first step and the relevant hypercohomology
    spectral sequence in the last.  Thus we have (a).

    Finally, if $T$ is concentrated in degree $0$, the exact triangle tells us that
    $H^i(G_{F,S},\F_{\Ga}(T))$ is finitely generated over $\F_{\Ga}(T)$ for sufficiently large $i$.
    Thus, we are still able to perform the inductive step in the proof of Theorem \ref{fg cohomology} 
    to obtain that $\R\Ga(G_{F,S},\F_{\Ga}(T))$ sits in $\DD^+_{\La-\ft}(\Mod_{\La})$
    for such $T$, and then for bounded below $T$ via the hypercohomology spectral
    sequence.  Thus, we have (b).
\epf

\br \label{fg coh local}
    Part (d) of Proposition \ref{fg type cohomology} makes sense and
    holds more generally for any complex $T$ in $\C_{R,G_v}^{R-\ft}$, as $\La_{\chi}$ has a 
    $G_v$-action through the composite with the canonical map $G_v \to G_{F,S}$.
\er

\br \label{nonarch case}
    We may also consider the setting in which we take $F$ itself to be a nonarchimedean local field
    of characteristic not equal to $p$.  We then let $G_F$ be its absolute Galois group,
    let $\Gamma$ denote the Galois group of
    a $p$-adic Lie extension of $F$, and as before, set $\La = R\ps{\Gamma}$.
    In this case, we obtain immediately from
    Theorem \ref{fg cohomology} that the cohomology groups $H^j(G_F,\F_{\Ga}(T))$ are finitely
    generated $\La$-modules for all $j$ and any complex of objects of
    $T$ of $\C_{R,G_F}^{R-\ft}$, and consequently,
    that $\R\Ga(G_F,\F_{\Ga}(T))$ is an object of $\DD_{\La-\ft}(\Mod_{\La})$ for any
    complex $T$ of objects in $\C_{R,G_F}^{R-\ft}$.
\er

\subsection{Duality over local fields}

In this subsection, we will state a version of Tate duality for $\Gamma$-induced modules
with $\La \ot_R^{\Li} \omega_R$ for a dualizing complex $\w_R$ replacing
$\Qp/\Zp$.  Fix a nonarchimedean prime $v \in S$.
We continue to 
suppose that $\Gamma$ is a compact $p$-adic Lie group that is a quotient
of $G_{F,S}$ defining an extension $F_{\infty}$, though as we remark later, 
we could just as
well assume it to be a quotient of $G_v$.
Let $T$ be a bounded complex of objects in $\C_{R,G_v}^{R-\ft}$.
As before, we let $J_R$ denote
a complex of injective $R$-modules that represents $\w_R$, and we let $T^*$ be a
quasi-isomorphic subcomplex of $\Hom_{R}(T,J_R)$ consisting of objects of $\C_{R,G_v}^{R-\ft}$.

The cup product of Lemma \ref{cupprod2} yields a map
$$
    C(G_v, \F_{\Ga}(T))\ot_{R} C(G_v, \F_{\Ga}(T^*)^{\iota}(1)) \lra
    \tau_{\ge 2}^{II}C(G_v, \La\ot_RJ_R (1))
$$
of $\La \ot_R \Lo$-modules.
Taking adjoints, we have the following maps of complexes of
$\La$-modules and $\Lo$-modules, respectively:
\begin{eqnarray*}
    &C(G_{v},\F_{\Ga}(T))\lra\Hom_{\Lo}\!\Big(C(G_{v},\F_{\Ga}(T^{*})^{\iota}(1)),\tau^{II}_{\geq
    2}C(G_{v}, \La\ot_R J_R (1))\Big),&\\
    &C(G_{v},\F_{\Ga}(T^*)^{\iota}(1))\lra\Hom_{\La}\!\Big(C(G_{v},\F_{\Ga}(T)),\tau^{II}_{\geq
    2}C(G_{v}, \La\ot_R J_R (1))\Big).&
\end{eqnarray*}

From now on, we use $M^{\vee}$ to denote the Pontryagin dual of a locally compact abelian group $M$.
We have the following analogue of \cite[(5.2.1)]{Ne}.

\bl \label{loctrunc}
    Let $A$ be a bounded complex in $\C_{\La-\La}^{\Lo-\ft}$ with terms that are flat as $R$-modules.
        Then there exists a quasi-isomorphism
        \[
                q_A \colon (A \ot_R J_R)[-2]  \lra \tau^{II}_{\geq 2}C(G_v, A \ot_R J_R (1)),
        \]
            natural in $A$,
        of chain complexes of $\La \ot_R \Lo$-modules.
\el

\begin{proof}
        We first note that if $M$ is any finitely generated $R$-module with a trivial action of $G_v$,
        then Tate local duality as in \cite[Theorem 4.1.2]{Lim-adic} yields
        isomorphisms
        $$
                H^2(G_v, A \ot_R M(1)) \cong H^0(G_v, (A \ot_R M)^{\vee})^{\vee}\cong A \ot_R M,
        $$
        of $\La \ot_R \Lo$-modules.  Now suppose that $M$ is any (ind-admissible)
        $R[G_v]$-module with trivial $G_v$-action, and write $M = \bigcup
        M_{\al}$, where the $M_{\al}$ are finitely generated $R$-submodules
        of $M$. Since the terms of $A$ are $R$-flat, $A^i \ot_R M$ is the union of the $A^i \ot_R M_{\al}$
        for any $i \in \Z$, so
        $$
                H^2(G_v, A \ot_R M(1)) \cong \ilim H^2(G_v, A \ot_R M_{\al}(1))\cong \ilim A \ot_R
                M_{\al} \cong A \ot_R M.
        $$
        As $G_v$ is of $p$-cohomological dimension $2$,
        we therefore have a quasi-isomorphism
        \[
            A \ot_R M[-2] \lra \tau_{\geq 2}^{II}C(G_v, A \ot_R M(1)),
        \]
        of complexes of $\La \ot_R \Lo$-modules.  The case of a bounded complex $M$ of
        ind-admissible $R$-modules, e.g. $M = J_R$, then follows easily.
        Naturality is immediate from the construction.
\end{proof}

Combining Lemma \ref{loctrunc} for $A = \La$
with the morphisms constructed above and passing to the derived
category, we obtain morphisms as in the following theorem.

\bt \label{local duality La}
    Let $T$ be a bounded complex in $\C^{R-\ft}_{R,G_v}$.  Then the morphisms
    \begin{eqnarray*}
    &\R\Ga(G_{v},\F_{\Ga}(T))\lra
    \RHom_{\Lo}\!\Big(\R\Ga(G_{v},\F_{\Ga}(T^{*})^{\iota}(1)),\La\ot^{\Li}_{R}\w_R\Big)[-2]&\\
    &\R\Ga(G_{v},\F_{\Ga}(T^*)^{\iota}(1))\lra
    \RHom_{\La}\!\Big(\R\Ga(G_{v},\F_{\Ga}(T)),\La\ot^{\Li}_{R}\w_R\Big)[-2]&
    \end{eqnarray*}
    in $\DD_{\La-\ft}(\Mod_{\La})$ and $\DD_{\La-\ft}(\Mod_{\Lo})$, respectively, are isomorphisms.
\et

Over the course of the next two subsections, we will show that the
first of the morphisms of Theorem \ref{local duality La} is an
isomorphism. The proof that the second is an isomorphism is
completely analogous.
By an argument similar to that in \cite[Lemma 5.2.5]{Ne}, it
suffices to prove the theorem for a particular translate of the
dualizing complex. In view of this, we will assume that our choice
of the dualizing complex satisfies \cite[(2.5)(i)]{Ne} throughout
the next two subsections.

\br In the setting that $F$ is a nonarchimedean local field of characteristic not equal to
$p$, that $F_{\infty}$ is any $p$-adic Lie extension of $F$, and
that $T$ is a bounded complex of objects in $\C_{R,G_F}^{R-\ft}$,
the argument we are about to describe
also carries over  as
in Remark \ref{nonarch case} to prove the direct analogue of Theorem
\ref{local duality La} with $G_v$ replaced by $G_F$. \er

\subsection{Change of rings} \label{change of rings}

Recall that $R$ is a complete commutative local Noetherian ring with
maximal ideal $\m$ and finite residue field $k$. In this subsection,
we shall show that it suffices to prove Theorem \ref{local duality
La} for any complete commutative local Noetherian ring $S$ that has
$R$ as a quotient.  (In this subsection only, $S$ denotes such a ring, rather than
the set of primes chosen above, which is used only implicitly.)

We suppose given a surjection $\psi \colon S \to R$, which necessarily induces an
isomorphism on residue
fields. We use $\m_S$ to denote the maximal ideal of $S$, and we set
$$
     d = \dim_k \mf{m}_S/\mf{m}_S^2 - \dim_k \mf{m}/\mf{m}^2.
$$
Note that if $A$ is any complex of $S$-modules, then $\Hom_S(R,A)$ is isomorphic to
the subcomplex of $A$ with terms the torsion submodules under $\ker \psi$.  In particular,
if $B$ is a complex of $R$-modules, then the resulting inclusion map induces an isomorphism
$$
    \Hom_R(B,\Hom_S(R,A)) \lri \Hom_S(B,A),
$$
as is easy enough to see by the usual adjointness principle.  In particular, we have
$$
    \Hom_R(B,R^{\vee}) \lri \Hom_S(B,S^{\vee}).
$$

With this comparison of Matlis duals with respect to $R$ and $S$ in hand, we now
compare Grothendieck duals with respect to these rings (see also \cite{LC}).
For this, we choose particular dualizing complexes $J_R$ and $J_S$ of injectives as follows.
Let $\bar{x}_1, \ldots, \bar{x}_r \in \mf{m}$ with images forming a $k$-basis of $\mf{m}/\mf{m}^2$.
Recall (cf.\ \cite[Section 0.4]{Ne}) that we may take $J_R$ to be the Matlis dual of
$C_R[r]$, where $C_R$ is defined to be the complex
$$
        \Bigg[
        R \lra \bigoplus_{i} R_{\bar{x}_i} \lra \bigoplus_{i<j} R_{\bar{x}_i\bar{x}_j}
        \lra \cdots \lra R_{\bar{x}_1 \cdots \bar{x}_r}
        \Bigg]
$$
in degrees $[0,r]$, with the usual \v{C}ech differentials, which is to say that
$$
    J_R = \Hom_R(C_R[r] ,R^{\vee}).
$$
We lift $\bar{x}_1, \ldots, \bar{x}_r$ to elements $x_1, \ldots, x_r$ of $\mf{m}_S$ and extend to
a sequence $x_1, \ldots, x_s$ such that the images of $x_1, \ldots, x_s$ form a basis of
$\mf{m}_S/\mf{m}_S^2$ and the elements $x_{r+1}, \ldots, x_s$ map trivially to $\mf{m}/\mf{m}^2$.
Modifying each $x_i$ with $i > r$ by an element of $\mf{m}_S^2$ if necessary, we may choose
these elements so that they in fact map trivially to $\mf{m}$.  We use this sequence to define $C_S$,
and we set $J_S = \Hom_S(C_S[s],S^{\vee})$.

\bl \label{base change}
    Let $T$ be a bounded complex in $\C^{R-\ft}_{R,G_v}$, which we view also
    as a complex in $\C^{S-\ft}_{S,G_v}$ via $\psi$.
    Then $\psi$ induces a natural map $J_R \to J_S[d]$ and in turn
    an isomorphism
    $$
         \Hom_R(T,J_R) \lri  \Hom_S(T,J_S)[d]
    $$
    of complexes in $\II_{S,G_v}$.
\el

\begin{proof}
    	The natural map $C_S \to C_R$ factors through
        	an isomorphism $R\ot_S C_S \xrightarrow{\sim} C_R$. This
        	induces an isomorphism
        	$$
                J_R \lri \Hom_S(R,J_S)[d],
        	$$
    	which is given by the following composition of isomorphisms
        	\begin{multline*}
            	\Hom_R(C_R[r], R^{\vee}) \lri \Hom_S(C_R[r], S^{\vee})\lri \Hom_S(R\ot_SC_S[r], S^{\vee}) \\
            	\lri \Hom_S(R, \Hom_S(C_S[r],S^{\vee})) \lri \Hom_S(R, \Hom_S(C_S[s],S^{\vee}))[d],
        	\end{multline*}
    	where the latter natural map is defined as in \cite[(1.2.15)]{Ne} (noting \cite[(1.2.5)]{Ne}).
        	The result then follows from the isomorphisms
        	$$
          	\Hom_R(T,J_R) \lri \Hom_R(T,\Hom_S(R,J_S))[d] \lri \Hom_S(T,J_S)[d].
        	$$
\end{proof}

We now proceed to the reductive step.

\bp \label{reduce to pow ser ring}
    Let $T$ be a bounded complex in $\C^{R-\ft}_{R,G_v}$.  Setting $\Om = S\ps{\Ga}$,
    there is a commutative diagram of natural morphisms
    $$
         \SelectTips{cm}{} \xymatrix{
        \R\Ga(G_v, \F_{\Ga}(T)) \ar[r] \ar[d]^{\wr} &
        \RHom_{\Lo}\!\Big(\R\Ga(G_{v},\F_{\Ga}(T^{*})^{\iota}(1)),\La\ot^{\Li}_{R}\w_R\Big)[-2]
        \ar[d]^{\wr} \\
        \R\Ga(G_v,\Om_{\chi} \ot_S T) \ar[r] &
        \RHom_{\Oo}\!\Big(\R\Ga(G_{v}, {}_{\chi}\Om \ot_S T^{*}(1)),\Om \ot_S^{\Li} \omega_S[d] \Big)[-2]
         }
    $$
    in which the horizontal morphisms are as in Theorem \ref{local duality La} and the
    vertical arrows are isomorphisms in $\DD(\Mod_{\Om})$.
    In particular, Theorem \ref{local duality La} holds for $R$ if it holds for $S$.
\ep

\begin{proof}
    Note that we have
    $$
        \Om_{\chi} \ot_S T \cong \La_{\chi} \ot_R T \quad (\mathrm{resp.,}\
    {}_{\chi} \Om \ot_S T^* \cong {}_{\chi} \La\ot_R
    T^*) $$
    so
    $\R\Ga(G_v, \F_{\Ga}(T))$ (resp., $\R\Ga(G_v, \F_{\Ga}(T^*)^{\iota}(1))$)
    is the same complex for $R$ and for $S$.
    By Lemma \ref{base change}, and employing the $S$-flatness of $\Om$ in
    applying Lemma \ref{tech 2},
    we have canonical isomorphisms
    $$
        \La \ot_R J_R \lri \Om \ot_S J_R
        \lri \Om \ot_S \Hom_S(R,J_S)[d] \lri \Hom_S(R,\Om
        \ot_S J_S)[d],
    $$
    of complexes of $\Oo \ot_S R \cong \Lo$-modules,
    and the composition of the composite isomorphism with the natural inclusion
    $$
        \Hom_S(R,\Om \ot_S J_S)[d] \lra \Om \ot_S J_S[d]
    $$
    induces the right-hand  vertical map in the proposition.  That the diagram commutes
    then follows directly from the definition of the cup product in Lemma \ref{cupprod2}.

    Let us fix a quasi-isomorphism $\iota \colon \Om \ot_S J_S \to K$, where
    $K$ is a bounded below complex of injective $\Oo$-modules.
    We check that the map
    $$
        \iota_R \colon \Hom_S(R,\Om \ot_S J_S) \lra \Hom_S(R,K)
    $$
    induced by $\iota$
    is a quasi-isomorphism, which implies that $\Hom_S(R, \Om \ot_S J_S)$ represents
    $\RHom_S(R,\Om \ot_S^{\Li} \omega_S)$ in $\DD(\Mod_{\Lo})$.
    Let $\e \colon P \to R$ be a resolution of $R$ by a complex $P$
    of finitely generated free $S$-modules.  We have a commutative
    diagram
    $$
        \SelectTips{cm}{} \xymatrix@C=20pt{
            \Om \ot_S \Hom_S(R,J_S) \ar[r]^-{\phi_R} \ar[d]^{\id \ot \e^*_J} &
            \Hom_S(R,\Om \ot_S J_S) \ar[r]^-{\iota_R} \ar[d]^{\e^*_{\Om \ot J}}
             & \Hom_S(R,K) \ar[d]^{\e^*_K} \ar[r]^-{\psi_R} &
            \Hom_{\Oo}(\Om \ot_S R, K) \ar[d]^{(\id \ot \e)^*_K} \\
            \Om \ot_S \Hom_S(P,J_S) \ar[r]^-{\phi_P} &
            \Hom_S(P,\Om \ot_S J_S) \ar[r]^-{\iota_P} & \Hom_S(P,K) \ar[r]^-{\psi_P} &
            \Hom_{\Oo}(\Om \ot_S P, K)
        }
    $$
    in which
    the vertical maps are all induced by $\e$.
    Moreover, note that the upper horizontal morphisms are of complexes of $\Lo$-modules
    and the others are of complexes of $\Oo$-modules.
    The maps $\phi_P$ and $\phi_R$ are isomorphisms, and
    the map $\id \ot \e_J^*$ is a quasi-isomorphism by the $S$-injectivity of $J_S$ and
    the $S$-flatness of $\Om$, so $\e^*_{\Om \ot J}$ is also a quasi-isomorphism.
    The maps $\psi_R$ and $\psi_P$ are isomorphisms by the usual adjointness of $\Hom$
    and the tensor product, and the map
    $(\id \ot \e)^*_K$ is a quasi-isomorphism by the $\Oo$-injectivity
    of $K$, so $\e_K^*$ is a quasi-isomorphism as well.
    Finally, the map $\iota_P$ is a quasi-isomorphism by the $S$-projectivity of $P$,
    and it follows that  $\iota_R$ is a quasi-isomorphism.

    To finish the proof, we need only show that the morphism
    $$
        	\RHom_{\Lo}\Big(X,\RHom_S(R,\Om \ot_S^{\Li} \w_S)\Big)
       	\lra \RHom_{\Oo}(X,\Om \ot_S^{\Li} \w_S)
    $$
    is an isomorphism for $X = C(G_v, \F_{\Ga}(T^*)^{\iota}(1))$    .
    It is easy to see that $\Hom_S(R,K)$
    is a complex of injective $\Lo$-modules.  Moreover, every $\Oo$-homomorphism from
    the complex $X$ of $\Lo$-modules to $K$ must factor through
    $$
        \Hom_{\Oo}(\La,K) \cong \Hom_{\Oo}(\Om \ot_S R,K) \cong \Hom_S(R,K)
    $$
    of $K$, so the map
    $$
        \Hom_{\Lo}(X,\Hom_S(R,K)) \lra \Hom_{\Oo}(X,K)
    $$
    induced by inclusion is an isomorphism, as desired.
\end{proof}

\subsection{Duality over flat $\Zp$-algebras}

In this subsection, we will prove Theorem \ref{local duality La}.
As remarked in the proof of Proposition \ref{Noetherian compact ring},
every complete commutative Noetherian local ring with
finite residue field is a quotient of a power series ring
in finitely many variables over an unramified extension of $\Zp$.
Proposition \ref{reduce to pow ser ring} then allows us to reduce
Theorem \ref{local duality La} to the case of such power series rings.
We therefore can and do assume that $R$ is $\Zp$-flat throughout this subsection.

We will require a bit more general of a derived adjoint
map than the one we wish to prove is an isomorphism. The
construction is found in the following easy lemma.

\bl \label{commutative derived}
    Let $A$ be a bounded complex of $R$-flat objects in $\C_{\La-\La}^{\Lo-\ft}$.  Then we have a morphism
    $$
        \R\Ga(G_v,A \ot_{\La} \F_{\Ga}(T)) \lra \RHom_{\Lo}\!\Big(\R\Ga(G_v,\F_{\Ga}(T^*)^{\iota}(1)),
        A \ot_R^{\Li} \w_R\Big)[-2]
    $$
    in $\DD(\Mod_{\La})$ that is natural in the complex $A$.
\el

\begin{proof}
    This is simply the derived adjoint of the composition of the cup product in Lemma \ref{cupprod2}
    with the map to the truncation $\tau^{II}_{\ge 2}C(G_v,A \ot_R J_R(1))$, which
    is quasi-isomorphic to $(A \ot_R J_R)[-2]$ by Lemma \ref{loctrunc}.
    Naturality is immediate.
\end{proof}

We will mostly be interested in complexes of length $2$, so we introduce
 the following notation.  Suppose that  $B$ is an object in $\C_{\La-\La}^{\Lo-\ft}$ and
that $A$ is a subobject of $B$, which
is to say that it is a closed $\La \cotimes{R} \Lo$-submodule with the subspace topology.
Then for any bounded complex
$T$ of objects in $\C_{R,G_v}^{R-\ft}$, we define a complex
\[ \F_{B/A}(T) = [A \to B]\ot_{\La}\F_{\Ga}(T), \]
of objects in $\C_{\La,G_v}$,
where $A$ and $B$ are in degree -1 and 0 respectively.

Let $I= I(\Ga)$ denote the augmentation ideal of $\La$.  We then have the following lemmas.
Recall for the first that we assume that $R$ is $\Zp$-flat in this subsection.

\bl \label{flat aug ideals} The following statements hold for any $n
\ge 0$.
\begin{enumerate}
    \item[$(a)$] The ideal $I^n$ is a flat $R$-module.
    \item[$(b)$] The module $I^n/I^{n+1}$ is finitely generated
    and of finite projective dimension over $R$.
\end{enumerate}
\el

\bpf
    (a) Suppose first that $\Ga$ is
    finite, and let $\mc{I}$ denote the augmentation ideal in $\Zp[\Ga]$.
    As $\Zp[\Ga]$ is finitely generated and free over $\Zp$, so is $\mc{I}^n$.
    Since $R$ is $\Zp$-flat, the natural surjection $R \ot_{\Zp} \mc{I}^n \to I^n$
    is an isomorphism, and therefore $I^n$ is free over $R$.

    In the general case,
    $\Ga \cong \plim \Ga'$, where $\Ga'$ runs over the Galois groups of the finite Galois
    extensions of $F$ in $F_{\infty}$, and
    $I^n$ is the inverse limit of the $n$th powers of the augmentation
    ideals $I(\Gamma')$ of the $R[\Ga']$.
    As $I(\Gamma')^n \cong I^n \ot_{\La} R[\Gamma']$,
    an application of Lemma \ref{completed tensor product} yields an isomorphism
    $$
        I^n \ot_R B \lri \plim_{\Gamma'}\,  ( I(\Gamma')^n \ot_R B)
    $$
    for any ideal $B$ of $R$.
    The left exactness of the inverse limit
    and the $R$-flatness of $I(\Ga')^n$ then imply that the canonical map
    $I^n \ot_R B \to I^n$
    is an injection, proving the flatness of $I^n$.

    \vspace{2ex}
    (b) Let $\mc{I}$ denote the augmentation ideal in $\Zp\ps{\Ga}$.  For each $n$, the composition
    \[ R \cotimes{\Zp}\mc{I}^n \lra R \cotimes{\Zp}\Zp\ps{\Ga} \lri R\ps{\Ga}\]
    induces a surjection
    $R \cotimes{\Zp}\mc{I}^n \tha I^n$ that fits into the
    following commutative diagram with exact rows:
    \[ \SelectTips{cm}{} \xymatrix@C=.4in{
      & R \cotimes{\Zp}\mc{I}^{n+1} \ar@{->>}[d] \ar[r]^{} & R \cotimes{\Zp} \mc{I}^n
      \ar@{->>}[d] \ar[r]^{} & R \cotimes{\Zp} (\mc{I}^n/\mc{I}^{n+1}) \ar[d]
      \ar[r] & 0 \\
      0 \ar[r] & I^{n+1} \ar[r] & I^{n} \ar[r] & I^{n}/I^{n+1} \ar[r] & 0.   }
    \]
    Since the two vertical maps on the left are surjections,
    so is the one on the right.
   Note that $\mc{I}^n/\mc{I}^{n+1}$ is a
    quotient of $n$th tensor power of the
    maximal abelian pro-$p$ quotient of $\Gamma$,
    a finitely generated $\Zp$-module. Therefore,
    the existence of the above surjection implies that
    $I^n/I^{n+1}$ is finitely generated over $R$. To see that the last
    assertion holds, we note that it follows from (a) that $I^n/I^{n+1}$ has finite flat
    dimension.
    Since flat dimension coincides with projective dimension for
    every finitely generated module over a Noetherian ring (see
    \cite[Proposition 4.1.5]{Wei}), we have our assertion.
\epf

\bl \label{coneisom}
    Let $T$ be a bounded complex of objects in $\C_{R,G_v}^{R-\ft}$.
    Then we have a quasi-isomorphism
    \[
        \Cone\Big(\F_{I^n/I^{n+1}}(T) \lra
        \F_{\La/I^{n+1}}(T)\Big) \lri \F_{\La/I^{n}}(T)
    \]
    of complexes in $\C_{\La,G_v}^{\La-\ft}$.
    Moreover, we have an exact triangle
    \[
        I^n/I^{n+1} \ot_R^{\Li} \omega_R \lra \La/I^{n+1} \ot_R^{\Li} \omega_R \lra
        \La/I^n \ot_R^{\Li} \omega_R
    \]
    in $\DD^b(\Mod_{\La \ot_R \Lo})$.
\el

\bpf
    Since the powers of $I$ are $R$-flat by Lemma \ref{flat aug ideals}(a),
    it suffices in both cases
    to show that there is a quasi-isomorphism
    \[ \Cone\big([ I^{n+1} \to I^n]
    \lra [ I^{n+1} \to \La]\big) \lri [I^n \to \La ].
    \]
    Note that the latter cone is precisely the complex
    \[ I^{n+1} \stackrel{f}{\lra} I^{n+1}\oplus I^n \stackrel{g}{\lra}
    \La,\] where $f(x) = (x,-x)$ and $g(x,y) = x+ y$ for $x\in I^{n+1}$
    and $y\in I^{n}$. One can now easily check that the diagram
    \[ \SelectTips{cm}{} \xymatrix @C=.72in{
        I^{n+1} \ar[d]_{} \ar[r]^-{f} & I^{n+1}\oplus I^n \ar[d]_{\nu} \ar[r]^-{g} & \La \ar@{=}[d] \\
        0 \ar[r]^{} & I^{n} \ar[r] & \La   }
    \]
    commutes, where $\nu$ is given by $\nu(x,y) = x+ y$ for $x\in
    I^{n+1}$ and $y\in I^{n}$, and the vertical maps induce isomorphisms
    on cohomology.
\epf

We are now able to prove the following proposition, which is an important ingredient in the proof of
Theorem \ref{local duality La}.

\bp \label{exact triangle}
    Let $T$ be a bounded complex of objects in $\C_{R,G_v}^{R-\ft}$. Then we have the following
    morphism of exact triangles
    \[ \SelectTips{cm}{} \xymatrix{
        \R\Ga(G_v, \F_{I^{n}/I^{n+1}}(T))\ar[d]_{} \ar[r]^{} & \RHom_{\Lo}\!\Big(\R\Ga(G_v,
        \F_{\Ga}(T^{*})^{\iota}(1)), I^{n}/I^{n+1}\ot_{R}^{\Li}\w_R \Big)[-2] \ar[d]^{} \\
        \R\Ga(G_v, \F_{\La/I^{n+1}}(T)) \ar[d]_{} \ar[r]^{} &  \RHom_{\Lo}\!\Big(\R\Ga(G_v,
        \F_{\Ga}(T^{*})^{\iota}(1)), \La/I^{n+1}\ot_{R}^{\Li}\w_R \Big)[-2] \ar[d]^{} \\
        \R\Ga(G_v, \F_{\La/I^{n}}(T)) \ar[r]^{} &  \RHom_{\Lo}\!\Big(\R\Ga(G_v,
        \F_{\Ga}(T^{*})^{\iota}(1)), \La/I^{n}\ot_{R}^{\Li}\w_R \Big)[-2]  }
    \]
    in $\DD(\Mod_{\La})$.
\ep

\bpf
    By Lemma
    \ref{flat aug ideals}(a), we see that
    $I^{n}/I^{n+1}\ot_{R}^{\Li}\w_R$ and $\La/I^{n}\ot_{R}^{\Li}\w_R$
    are represented by $[ I^{n+1} \to I^n]\ot_RJ_R $ and
    $[I^n \to \La]\ot_RJ_R $ respectively. Therefore, the
    commutativity of the diagram in the proposition follows from the naturality in Lemma
    \ref{commutative derived}.
    By Lemma \ref{coneisom}, both columns are exact triangles.
\epf

We now describe the idea of the proof of Theorem \ref{local duality
La}. We shall first prove that the morphism
$$
    \R\Ga(G_v, \F_{\La/I^n}(T)) \lra \RHom_{\Lo}\!\Big(\R\Ga(G_v,
    \F_{\Ga}(T^{*})^{\iota}(1)), \La/I^n\ot_{R}^{\Li}\w_R \Big)[-2]
$$
is an isomorphism for all $n$. Then, Theorem \ref{local duality La}
will follow from this by a limit argument. To show that above
morphism is an isomorphism, we will utilize Proposition \ref{exact
triangle}. Note that if any two of the horizontal morphisms in
Proposition \ref{exact triangle} is a quasi-isomorphism, so is the
third one. Therefore, by an inductive argument, we are reduced to
showing that the morphism
$$
    \R\Ga(G_v, \F_{I^{n}/I^{n+1}}(T)) \lra \RHom_{\Lo}\!\Big(\R\Ga(G_v,
    \F_{\Ga}(T^{*})^{\iota}(1)), I^{n}/I^{n+1}\ot_{R}^{\Li}\w_R\Big)[-2]
$$
is an isomorphism for all $n\geq 0$.

Note that $\Ga$
acts trivially on $I^{n}/I^{n+1}$. Therefore, one may view
$I^{n}/I^{n+1}$ as a $\Lo$-module via the augmentation map
$\La\to R$. We now have the following lemma.

\bl \label{duality lemma}
    Let $T$ be a bounded complex of objects in
    $\C_{R,G_v}^{R-\ft}$. Then we have the following isomorphisms
    \[
        \F_{I^n/I^{n+1}}(T) \stackrel{\sim}{\longleftarrow}
        I^n/I^{n+1} \dotimes{\La} \F_{\Ga}(T) \lri
        I^n/I^{n+1} \dotimes{R} T
    \]
    in $\DD^b(\C_{\La,G_v})$. Therefore,  we have an isomorphism
    \[
        \R\Ga(G_v, \F_{I^n/I^{n+1}}(T))\lri
        \R\Ga(G_v, I^n/I^{n+1}
    \dotimes{R}
        T)
    \]
    in $\DD(\Mod_\La)$.
\el

\bpf
    Let $P$ be a resolution of $I^n/I^{n+1}$ consisting of
    projective objects in $\C_{\La-\La}$.
    Then there is a quasi-isomorphism $P \to [ I^{n+1} \to I^{n} ]$
    of complexes in $\C_{\La-\La}$ which lifts the identity map on $I^n/I^{n+1}$.
    By Lemmas \ref{projobj} and \ref{flat aug ideals}(a) this
    induces a quasi-isomorphism
        \[
            P\ot_{\La}\La_{\chi} \ot_{R}T \lra [ I^{n+1} \to I^{n} ]\ot_{\La}\La_{\chi}\ot_{R} T
    \]
    of bounded above complexes in $\C_{\La,G_v}$, proving the first isomorphism.
    The second isomorphism is a special case of Lemma \ref{bbc} with $N = I^n/I^{n+1}$
    and $\Gamma'$ trivial.
    Note that both quasi-isomorphisms of complexes factor
    through a quasi-isomorphic truncation of $P\ot_{\La}\La_{\chi} \ot_{R}T$, so they may be seen
    in the bounded derived category.
\epf

\bl \label{tech diagram}
    For each $n$, there is a commutative diagram
    \[\SelectTips{cm}{} \xymatrix {
        I^n/I^{n+1}\ot_R^{\Li} \R\Ga(G_v,T) \ar[d]_{\wr} \ar[r] &
        I^n/I^{n+1}\ot_R^{\Li} \RHom_R\!\Big(\R\Ga(G_v,T^*(1)), \w_R\Big)[-2] \ar[d]^{\wr} \\
        \R\Ga(G_v, I^n/I^{n+1}\dotimes{R} T) \ar[dd]_{\wr} \ar[r]^{} &
        \RHom_R\!\Big(\R\Ga(G_v,T^*(1)), I^n/I^{n+1}\ot_R^{\Li}\w_R\Big)[-2] \ar[d]^{\wr} \\
                  & \RHom_R\!\Big(\R\Ga(G_v,\F_{\Ga}(T^*)^{\iota}(1))\ot_{\La}^{\Li}R,
        I^{n}/I^{n+1}\ot_R^{\Li}\w_R\Big)[-2] \ar[d]^{\wr} \\
        \R\Ga(G_v, \F_{I^{n}/I^{n+1}}(T)) \ar[r]
        & \RHom_{\Lo}\!\Big(\R\Ga(G_v,\F_{\Ga}(T^*)^{\iota}(1)), I^n/I^{n+1}\ot_R^{\Li}\w_R\Big)[-2],
         }
    \]
    where the vertical morphisms are isomorphisms in $\DD(\Mod_R)$.
\el

\bpf
        By Lemma \ref{flat aug ideals}(b), we may choose a bounded resolution $Q$ of
    $I^n/I^{n+1}$ by finitely generated projective $R$-modules.
        By \cite[Proposition 3.4.4]{Ne}, we have an isomorphism of complexes
        \[
                \al \colon Q\ot_R C(G_v, T) \lri C(G_v, Q\ot_R T)
        \]
        that fits into the commutative diagram
        \[
        \SelectTips{cm}{} \xymatrix{
            Q\ot_R C(G_v, T)\ot_R C(G_v, T^*(1)) \ar[d]_{\al\ot\id} \ar[r] &
        Q\ot_R C(G_v, J_R(1)) \ar[d]^{\al'} \\
            C(G_v, Q\ot_R T)\ot_R C(G_v, T^*(1)) \ar[r] & C(G_v, Q\ot_R J_R(1)),
            }
    \]
    where $\al'$ is defined analogously to $\alpha$ and is also an isomorphism.
    Since $J_R$ is a bounded complex of $R$-injectives,
    we may find homotopy inverses to the quasi-isomorphisms
    $q_R$ and $q_Q$ of Lemma \ref{loctrunc} that fit in a commutative diagram
        \[
        \SelectTips{cm}{} \xymatrix{
            Q\ot_R C(G_v, J_R(1)) \ar[d]^{\al'} \ar[r]^{} & Q\ot_R \tau^{II}_{\geq 2}C(G_v, J_R(1))
            \ar[d]_{} \ar[r] &  Q\ot_RJ_R  [-2] \ar@{=}[d]\\
            C(G_v, Q\ot_RJ_R(1))\ar[r] &  \tau^{II}_{\geq 2}C(G_v, Q\ot_RJ_R(1)) \ar[r]
            & Q\ot_RJ_R [-2]. }
        \]
        We obtain the top commutative square in the lemma and the fact that the
        vertical morphisms therein are isomorphisms, the one on the right by Lemma
        \ref{tech 2 derived} and \cite[Proposition 4.2.3]{Ne}.

        Let $P$ be a resolution of $I^n/I^{n+1}$ consisting of finitely
        generated projective $\Lo$-modules.
        We may view $Q$ as a resolution of $\Lo$-modules via the augmentation map
        $\La\to R$.  Then, by Proposition \ref{derived comp tensor}, the map
        \[
                P\ot_{\La}\F_{\Ga}(T) =
                P\ot_{\La}\La_{\chi} \ot_{R}T \lra Q\ot_{\La}\La_{\chi} \ot_{R}T
                \cong Q\ot_{R}T
        \]
        is a quasi-isomorphism of complexes of objects in $\C_{R,G_v}$,
            and we let $f$ denote its induced map on cochains.
        Let $L$ be a resolution of $R$ consisting of
        finitely generated projective $\La$-modules.
        Then, by an opposite version of Theorem \ref{descseq},
        we have a quasi-isomorphism
        \[
                g \colon C(G_v,\F_{\Ga}(T^*)^{\iota}(1))\ot_{\La}L  \lra  C(G_v,T^*(1))
        \]
        of complexes of $R$-modules.  Noting Lemma \ref{cupprod2}, we have a commutative diagram
    \[
        \small
        \SelectTips{cm}{} \xymatrix@C=18pt@R=18pt{
                C(G_v, P\ot_{\La}\F_{\Ga}(T))\ot_R C(G_v, \F_{\Ga}(T^*)^{\iota}(1)) \ot_{\La}L
                \ar[d]_{f \ot g} \ar[r]^{} &  C(G_v, P \ot_R J_R (1))\ot_{\La}L\ar[d]^{\e} \\
                C(G_v, Q\ot_{R}T)\ot_R C(G_v, T^*(1)) \ar[r]^{} & C(G_v, Q\ot_RJ_R (1)),   }
        \]
        where $\e$ is induced by the augmentation $L\to R$
        and the map $P \to Q$.

        Taking adjoints and applying the homotopy inverse to $q_Q$
    as above, we obtain  the commutative diagram
        \[\small
        \SelectTips{cm}{} \xymatrix@C=18pt@R=18pt{
                C(G_v, Q\ot_R T)  \ar[r]^{} & \Hom_R\!\Big(C(G_v,T^*(1)),
        Q \ot_R J_R  \Big)[-2] \ar[d] \\
                C(G_v, P\ot_{\La}\F_{\Ga}(T)) \ar[u] \ar @{=}[d] \ar[r]^{} &
                \Hom_R\!\Big( C(G_v, \F_{\Ga}(T^*)^{\iota}(1)) \ot_{\La}L,
                Q \ot_R J_R \Big)[-2] \ar[d]^{\wr}  \\
                C(G_v, P\ot_{\La}\F_{\Ga}(T)) \ar[r]^{} &
                \Hom_{\Lo}\!\Big( C(G_v, \F_{\Ga}(T^*)^{\iota}(1)),
                \Hom_R(L,Q \ot_R J_R )\Big)[-2],}
        \]
        which yields upon passage to the derived category
        the lower part of the diagram in the statement of the lemma,
        the vertical morphisms therein being isomorphisms by Theorem \ref{descseq},
        Lemma \ref{duality lemma} and Lemma \ref{derived adj}.
\epf

\bl
    The morphisms
    $$
        \R\Ga(G_v, \F_{I^{n}/I^{n+1}}(T)) \lra \RHom_{\Lo}\!\Big(\R\Ga(G_v,
        \F_{\Ga}(T^{*})^{\iota}(1)), I^{n}/I^{n+1}\ot_{R}^{\Li}\w_R\Big)[-2]
    $$
    in $\DD(\Mod_{\La})$ are isomorphisms for every $n\geq 0$.
\el

\bpf
   By Corollary \ref{derived Hom res}, it suffices to show that the
   above morphism is an isomorphism in $\DD(\Mod_R)$.
    The morphism
    \[
        \R\Ga(G_v,T)\lra \RHom_R\!\Big(\R\Ga(G_v,T^*(1)),\w_R\Big)[-2]
    \]
    in $\DD(\Mod_R)$ is an isomorphism by \cite[Proposition 5.2.4(ii)]{Ne},
    and so the top morphism of the diagram in Lemma \ref{tech diagram}
    is an isomorphism. Since all the vertical morphisms in the diagram
    are isomorphisms, it follows that the bottom morphism is also an
    isomorphism, as required.
\epf

\bp \label{n step}
    The morphisms
    $$
        \R\Ga(G_v, \F_{\La/I^{n}}(T)) \lra \RHom_{\Lo}\!\Big(\R\Ga(G_v,
        \F_{\Ga}(T^{*})^{\iota}(1)), \La/I^{n}\ot_{R}^{\Li}\w_R\Big)[-2]
    $$
    in $\DD(\Mod_{\La})$ are isomorphisms for every $n\geq 1$.
\ep

\bpf As seen in the above discussion, the preceding lemma allows
us to perform an inductive argument using the morphism of exact
triangles in Proposition \ref{exact triangle} to obtain the required
conclusion.
\epf

We finish the proof of Theorem \ref{local duality La} by passing
to the inverse limit.

\begin{proof}[Proof of Theorem \ref{local duality La}]
    By Remark \ref{fg coh local}, there exists a quasi-isomorphism
    $$
        W \lra C(G_v,\F_{\Ga}(T^{*})^{\iota}(1)).
    $$
    with $W$ a bounded above complex of finitely generated projective $\Lo$-modules.
    Since $J_R$ has cohomology groups which are finitely generated over $R$,
    Lemma \ref{complex finite cohom} implies the existence of
    a subcomplex $C$ of $J_R$ such that $C$ is
    a complex of finitely generated $R$-modules and the inclusion $i:
    C \hra J_R$ is a quasi-isomorphism. We fix such a $C$ and write $X_n$ for $[I^{n} \to \La]$.
    The complex $\Hom_{\Lo}(W, X_n \ot_R C)$ represents
    \[
        \RHom_{\Lo}\!\Big(\R\Ga(G_v, \F_{\Ga}(T^{*})^{\iota}(1)),
        \La/I^{n}\ot_{R}^{\Li}\w_R \Big),
    \]
    and $\Hom_{\Lo}(W, \La \ot_R C)$ represents
    \[
        \RHom_{\Lo}\!\Big(\R\Ga(G_v, \F_{\Ga}(T^{*})^{\iota}(1)), \La\ot_{R}^{\Li}\w_R \Big),
    \]
    since $X_n$ is a complex of flat $R$-modules by Lemma \ref{flat aug ideals}(a).
    Now, for each $n$, we have a commutative diagram
    \[ \SelectTips{cm}{} \xymatrix{
        \R\Ga(G_v, \F_{\Ga}(T)) \ar[d]_{} \ar[r]^{} & \RHom_{\Lo}\!\Big(\R\Ga(G_v,
        \F_{\Ga}(T^{*})^{\iota}(1)), \La\ot_{R}^{\Li}\w_R\Big) [-2] \ar[d]^{} \\
        \R\Ga(G_v, \F_{\La/I^{n}}(T)) \ar[r]^{} & \RHom_{\Lo}\!\Big(\R\Ga(G_v,
        \F_{\Ga}(T^{*})^{\iota}(1)), \La/I^{n}\ot_{R}^{\Li}\w_R\Big)[-2]    }
    \]
    which induces a commutative diagram
    \[ \SelectTips{cm}{} \xymatrix{
        H^i(G_v, \F_{\Ga}(T)) \ar[d]_{} \ar[r]^{} &  H^{2-i}(\Hom_{\Lo}(W, \La \ot_R C)) \ar[d]^{} \\
        H^i(G_v, \F_{\La/I^{n}}(T)) \ar[r]^{} &  H^{2-i}(\Hom_{\Lo}(W, X_n \ot_R C))  }
    \]
    of cohomology groups. Since the maps in this diagram are compatible as we vary $n$,
    we obtain the commutative diagram
    \[ \entrymodifiers={!! <0pt, .7ex>+} \SelectTips{cm}{} \xymatrix{
        H^i(G_v, \F_{\Ga}(T)) \ar[d]_{} \ar[r]^{} &  H^{2-i}(\Hom_{\Lo}(W, \La \ot_R C)) \ar[d]^{} \\
        \plim_n H^i(G_v, \F_{\La/I^{n}}(T)) \ar[r]^{} & \plim_n H^{2-i}(\Hom_{\Lo}(W, X_n \ot_R C)).  }
    \]

    It remains to show that the upper horizontal map in the latter diagram
    is an isomorphism. By
    Proposition \ref{n step}, we know that the lower horizonal map is one.
    Noting that the maps $X_{n+1} \ot_R C \to X_n \ot_R C$ are injections of complexes
    with intersection $\La \ot_R C$ over all $n$,
    we have an isomorphism
    \[
        \plim_n\Hom_{\Lo}(W, X_n \ot_R C)
        \cong \Hom_{\Lo}(W, \La \ot_R C)
    \]
    of complexes of finitely generated $\La$-modules.
    Since inverse limits are exact for finitely generated $\La$-modules
    (being that they are compact),
        after taking cohomology groups, we have that
     the vertical map on the right is an
    isomorphism. On the other hand, we also have an isomorphism
    \[
        \plim_n \, (X_n\ot_{\La} \F_{\Ga}(T)) \cong \F_{\Ga}(T)
    \]
    of complexes of objects in $\C^{\La-\ft}_{\La,G_v}$ and hence an
    isomorphism
    \[
        \plim_n C(G_{v}, X_n\ot_{\La} \F_{\Ga}(T))\cong C(G_{v},\F_{\Ga}(T)).
    \]
    By \cite[Proposition 3.2.13]{Lim-adic},
    we then have
    \[
        \plim_n H^{j}(G_{v}, X_n\ot_{\La} \F_{\Ga}(T))\cong H^{j}(G_{v},\F_{\Ga}(T)).
    \]
    Therefore, the vertical map on the left is also an isomorphism. Hence, the top map is an
    isomorphism, as required.
\epf

\subsection{Duality over global fields} \label{global duality}

We end this paper by describing the global analog of Theorem
\ref{local duality La} that is our main result. Let $T$ be a bounded
complex of objects in $\C^{R-\ft}_{R,G_{F,S}}$, and we choose $T^*$
in $\C^{R-\ft}_{R,G_{F,S}}$ as in Section \ref{duals}. As in
\cite[(5.3.3)]{Ne}, we define two morphisms of complexes of $\La
\ot_R \Lo$-modules
\begin{eqnarray*}
    &_{c}\cup \colon C_{(c)}(G_{F,S},\F_{\Ga}(T))\ot_{R}C(G_{F,S}, \F_{\Ga}(T^{*})^{\iota}(1))\lra
    C_{(c)}(G_{F,S},\La\ot_{R}J_R (1))&\\
    &\cup_{c} \colon C(G_{F,S},\F_{\Ga}(T))\ot_{R}C_{(c)}(G_{F,S}, \F_{\Ga}(T^{*})^{\iota}(1))\lra
    C_{(c)}(G_{F,S},\La\ot_{R}J_R (1))&
\end{eqnarray*}
which are given by the formulas
\begin{eqnarray*}
    &(a,a_{S})\, {}_{c}\!\cup b = (a\cup b, a_{S}\cup_S \res_{S}(b))& \\
    &a\cup_{c}(b,b_{S}) = (a\cup b,(-1)^{\bar{a}}\res_{S}(a)\cup_S b_{S}),&
\end{eqnarray*}
where $\bar{a}$ denotes the degree of $a$, the direct sum of the restrictions
to the primes in $S$ is denoted $\mr{res}_S$, the symbol $\cup$ is the total cup product
\[
    C(G_{F,S},\F_{\Ga}(T))\ot_{R}C(G_{F,S}, \F_{\Ga}(T^{*})^{\iota}(1))\lra
    C(G_{F,S},\La\ot_{R}J_R (1))
\]
of Lemma \ref{cupprod2}, and $\cup_{S}$ is the direct sum of the corresponding local cup products.
Much as in Lemma \ref{loctrunc} (but now analogously to \cite[Lemma 5.7.3]{Ne}), we have a
quasi-isomorphism
\[
    \La \ot_R J_R  \lra \tau^{II}_{\geqslant 3}C_{(c)}(G_{F,S},\La\ot_RJ_R (1))
\]
of complexes of $\La \ot_R \Lo$-modules, which allows us to use adjoints to define the morphisms in the following theorem.

\bt \label{Global duality La}
    Let $T$ be a bounded complex of objects in $\C_{R,G_{F,S}}^{R-\ft}$.  Then
        we have an isomorphism
    $$
        \SelectTips{cm}{}  \entrymodifiers={!! <0pt, .8ex>+}  \xymatrix{
                \R\Ga_{(c)}(G_{F,S},\F_{\Ga}(T)) \ar[r]^-{\sim} \ar[d] &
        \RHom_{\Lo}\!\Big(\R\Ga(G_{F,S},\F_{\Ga}(T^*)^{\iota}(1)),\La \ot_R^{\Li} \w_R\Big)[-3] \ar[d] \\
                \R\Ga(G_{F,S},\F_{\Ga}(T)) \ar[r]^-{\sim} \ar[d] &
        \RHom_{\Lo}\!\Big(\R\Ga_{(c)}(G_{F,S},\F_{\Ga}(T^*)^{\iota}(1)),\La \ot_R^{\Li} \w_R\Big)[-3] \ar[d] \\
                \displaystyle\bigoplus_{v \in S} \R\Ga(G_v,\F_{\Ga}(T)) \ar[r]^-{\sim}  &
              \displaystyle\bigoplus_{v \in S} \RHom_{\Lo}\!\Big(\R\Ga(G_v,\F_{\Ga}(T^*)^{\iota}(1)),
              \La \ot_R^{\Li} \w_R\Big)[-2]
        }
    $$
     of exact triangles in $\DD_{\La-\ft}(\Mod_{\La})$.
\et

\bpf
        The morphism of exact triangles can be constructed as in
        \cite[Theorem 4.2.6]{Lim-adic}.  We remark that it suffices to prove that any two of
    the three horizontal morphisms is an isomorphism,
    and we focus on the second and third.
    By an analogous result to Proposition \ref{reduce to pow ser ring},
    we may assume that $R$ is regular and $\Zp$-flat.
    That the second of the morphisms
    is an isomorphism follows as in proof of Theorem \ref{local duality La}.
    In fact, aside from the obvious changes of notation for cochains and other objects, there
    are no significant changes to the proof.
    As for the third, note that for a nonarchimedean prime $v$ in $S$, the map
        $$
                \tau_{\ge 2}^{II} C(G_v,\Lambda \otimes_R J_R(1))
                \lra \tau_{\ge 2}^{II} \Big(C_{(c)}(G_{F,S},\Lambda \otimes_R J_R(1))[1]\Big)
        $$
        is a quasi-isomorphism that induces the identity maps on the cohomology groups of each complex.
        Under this identification, the map on the $v$-summand in the lower horizontal
        morphism in the statement is exactly the
        first isomorphism of Theorem \ref{local duality La}.

    It remains to prove that, for a real place $v$, the morphisms
        $$
        	\R\Ga(G_{v},\F_{\Ga}(T))\lra
                \RHom_{\Lo}\!\Big(\R\Ga(G_{v},\F_{\Ga}(T^{*})^{\iota}(1)),\La\ot^{\Li}_{R}\w_R\Big)[-2]
        $$
        are isomorphisms in $\DD_{\La-\ft}(\Mod_{\La})$.  There are two cases:
    either $v$ extends to a complex place or splits completely in $F_{\infty}$.
    However, if $v$ becomes complex,
    then both $\widehat{C}(G_v,\F_{\Ga}(T))$ and
    $\widehat{C}(G_v,\F_{\Ga}(T^{*})^{\iota}(1))$ are acyclic, and there is nothing to prove.

    If $v$ splits completely in $F_{\infty}$, then
    the adjoint map of interest is identified (as in the proof of Proposition \ref{spectri})
    with
    \begin{eqnarray*}
        &\La \otimes_R \widehat{C}(G_v,T) \lra
        \Hom_{\Lo}\!\Big(\La \ot_R \widehat{C}(G_v,T^*(1)), \La \ot_R J_R\Big)[-2],&\\
        &\la \ot f \mapsto (\mu \otimes g \mapsto \la\mu \otimes f \cup g),&
    \end{eqnarray*}
    where
    $$
        \cup \colon \widehat{C}(G_v,T) \otimes_R \widehat{C}(G_v,T^*(1)) \lra J_R[-2]
    $$
    is given by cup product followed by the maps
    $$
        \widehat{C}(G_v,J_R(1)) \lra \tau_{\ge 2}^{II}\widehat{C}(G_v,J_R(1)) \lra
        \tau_{\ge 2}^{II} \Big(C_{(c)}(G_{F,S},J_R(1))[1]\Big)
        \lra J_R[-2],
    $$
    the latter being a homotopy inverse to the inclusion of complexes.

    Suppose now that we have a quasi-isomorphism
    $D \to \widehat{C}(G_v, T^*(1))$ with $D$ a $q$-projective
    complex of finitely generated $R$-modules.
    In the commutative diagram
    $$
        \SelectTips{cm}{} \xymatrix@C=20pt{
         \Hom_{\Lo}\!\Big(\La \ot_R \widehat{C}(G_v,T^*(1)),
        \La \ot_R J_R\Big) \ar[d]^{\wr} \ar[r] & \Hom_{\Lo}(\La \ot_R D,\La \ot_R J_R) \ar[d]^{\wr} \\
        \Hom_R\!\Big(\widehat{C}(G_v,T^*(1)), \La \ot_R J_R\Big) \ar[d] \ar[r] &
        \Hom_R(D,\La \ot_R J_R)
        \ar[d]^{\wr}  \\
        \La \ot_R \Hom_R\!\Big(\widehat{C}(G_v,T^*(1)), J_R\Big) \ar[r]  & \La \ot_R \Hom_R(D,J_R),
        }
    $$
    of canonical maps,
    the lower right vertical map is an isomorphism by Lemma \ref{tech 2}
    (the terms of $D$ being finitely generated over $R$), and
    the lower horizontal morphism is a quasi-isomorphism by the $R$-injectivity of the terms of $J_R$
    and the $R$-flatness of $\La$.  It follows that the derived adjoint
    morphism is represented by a morphism
    $$
        \La \otimes_R \widehat{C}(G_v,T) \lra \La \ot_R \Hom_R\!\Big(\widehat{C}(G_v,T^*(1)),J_R\Big)[-2]
    $$
    that by construction takes $\la \otimes f$ to $\la \otimes (g \mapsto f \cup g)$, which
    is the tensor product of the identity map with the derived adjoint that is already known
    to be a quasi-isomorphism by \cite[(5.7.5)]{Ne}.

    We are left to construct $D$, the existence of which is not so obvious as
    $\widehat{C}(G_v, T^*(1))$ need not have bounded cohomology.  However, as is seen
    in \cite[Lemma 3.3]{Sp}
    (see also \cite[Appendix]{Kel}),
    such a complex
    may be constructed as a direct limit $\ilim D_n$ using
    maps $\alpha_n \colon D_n \to D_{n+1}$
    that are split injective in each degree, and
    where
        $$
        D_0 \lra \tau_{\le 0}^{II} \widehat{C}(G_v,T^*(1))
    $$
    is a quasi-isomorphism from
    a bounded above complex of finitely generated $R$-projectives and
    $E_n = D_n/\alpha_{n-1}(D_{n-1})$ for $n \ge 1$
    is chosen to be quasi-isomorphic to
    $$
        \Cone \Big( \tau_{\le n-1}^{II} \widehat{C}(G_v, T^*(1)) \lra \tau_{\le n}^{II} \widehat{C}(G_v,T^*(1))
        \Big),
    $$
    which has bounded, $R$-finitely generated cohomology.
    Since $R$ has finite global dimension, we may choose the $E_n$ as bounded complexes
    of finitely generated
    projective $R$-modules such that, for each $k$, only finitely many $E_n^k$ are nonzero,
    which is what was needed.
\epf

We end with a straightforward remark.

\begin{remark}
    There is also an analogous diagram and isomorphisms
    for the other adjoint, with morphisms as in the
    second map in Theorem \ref{local duality La}, and this follows, for instance,
    by a nearly identical argument.
\end{remark}


\end{document}